\documentclass[11pt,a4paper]{article}
\usepackage[T1]{fontenc}
\usepackage{lmodern}
\usepackage{amsmath,amssymb,amsthm,mathtools}
\usepackage[margin=28mm,headheight=14pt]{geometry}
\usepackage{array,longtable,booktabs}
\usepackage{xcolor}
\usepackage[colorlinks=true,linkcolor=blue!50!black,citecolor=blue!50!black,urlcolor=blue!50!black]{hyperref}
\usepackage{enumitem}
\setlist{itemsep=3pt,topsep=5pt}
\numberwithin{equation}{section}
\newtheorem{theorem}{Theorem}[section]
\newtheorem{lemma}[theorem]{Lemma}
\newtheorem{proposition}[theorem]{Proposition}
\newtheorem{corollary}[theorem]{Corollary}
\theoremstyle{definition}
\newtheorem{definition}[theorem]{Definition}

\theoremstyle{remark}
\newtheorem{remark}[theorem]{Remark}
\newcommand{\R}{\mathbb R}

\newcommand{\Sym}{\operatorname{Sym}}
\newcommand{\End}{\operatorname{End}}
\newcommand{\tr}{\operatorname{tr}}
\newcommand{\im}{\operatorname{im}}
\newcommand{\rank}{\operatorname{rank}}
\newcommand{\Ann}{\operatorname{Ann}}
\newcommand{\Isom}{\operatorname{Isom}}
\newcommand{\Kill}{\operatorname{Kill}}
\newcommand{\so}{\mathfrak{so}}
\newcommand{\cyc}{\operatorname{cyc}}
\newcommand{\Id}{\operatorname{Id}}

\newcommand{\Lie}{\mathcal L}
\newcommand{\dd}{\mathrm d}
\newcommand{\ip}[2]{\langle #1,#2\rangle}
\newcommand{\norm}[1]{\lVert #1\rVert}
\newcommand{\nablah}{\nabla^h}
\newcommand{\nablac}{\nabla^{\mathrm c}}
\newcommand{\Cent}{\operatorname{Cent}}
\newcommand{\rad}{\operatorname{rad}}
\hypersetup{pdftitle={Rigidity of Closed Landsberg Manifolds in Higher Dimensions},pdfauthor={Zhuoqi Li, Jianyu Mao, Linfeng Zhou}}
\title{Rigidity of Closed Landsberg Manifolds\\in Higher Dimensions}
\author{Zhuoqi Li \and Jianyu Mao \and Linfeng Zhou}
\date{}
\begin{document}
\maketitle
\begin{abstract}
We prove that a regular Landsberg metric on a closed manifold of dimension at least four is Berwald. The proof uses the intrinsic isometry group of the indicatrix. In the nontransitive case, the geometry of its orbits and the vertical Bianchi identity imply that the nonlinear curvature is linear in the fibre variable. In the transitive case, a rigidity theorem for the angular Berwald quartic reduces the problem to the same linearity statement on an indicatrix of constant sectional curvature one. The metric is therefore $R$-quadratic, and Crampin's compact rigidity theorem applies. The angular quartic theorem and the required fibrewise linearity results are proved in the appendices.
\end{abstract}

\section{Introduction}\label{sec:introduction}

The Landsberg--Berwald problem asks whether every regular Landsberg metric is Berwald. For a Landsberg metric, the differential of nonlinear parallel translation preserves the fundamental tensor. For a Berwald metric, parallel translation is induced by a linear map between tangent spaces. The implication from the second condition to the first is classical; the converse is the rigidity problem considered here.

Throughout this paper, \emph{closed} means compact without boundary. A Finsler metric is called \emph{regular} when it is smooth and strongly convex in every nonzero tangent direction. We prove the following theorem.

\begin{theorem}\label{thm:main}
Let $M$ be a closed smooth manifold of dimension $n\geq4$, and let $F$ be a regular Landsberg metric on $M$. Then $F$ is Berwald.
\end{theorem}

The precise regularity conventions are recalled in Section~\ref{sec:preliminaries}. Neither $F(x,-y)=F(x,y)$ nor real analyticity in $y$ is assumed.

For the standard definitions and identities in Finsler geometry, we refer to Bao--Chern--Shen~\cite{BCS}. Earlier rigidity theorems require additional curvature or structural hypotheses. These include Numata's theorem for scalar curvature, Matsumoto's subsequent results, and the work of B\'acs\'o--Ilosvay--Kis on the Douglas condition~\cite{Numata,Matsumoto,Bacso}. Rigidity has also been established for special metric classes, notably regular $(\alpha,\beta)$-metrics and $(\alpha_1,\alpha_2)$-metrics~\cite{ShenAB,XuDeng}. Xu--Matveev proved rigidity for tangent norms with linear $SO(k)\times SO(n-k)$ symmetry~\cite{XuMatveev}. Their hypothesis prescribes linear norm isometries, whereas the group used below consists of intrinsic isometries of the indicatrix metric.

Shen introduced the class of $R$-quadratic metrics~\cite{ShenR}. Crampin proved the compact $R$-quadratic Landsberg rigidity theorem~\cite{Crampin2011}; see also the survey~\cite{Crampin2024}. This theorem supplies the final step of our argument in Section~\ref{sec:compact}.

Mao--Zhou~\cite{MaoZhou} proved rigidity for closed three-dimensional regular Landsberg manifolds. Their argument separates indicatrices according to the dimension of the Killing algebra; we instead distinguish the transitive and nontransitive cases of the intrinsic indicatrix action. This gives a common framework in dimensions at least three, with the angular quartic theorem providing the additional ingredient needed in higher dimensions. Li--Liu~\cite{LiLiu} established rigidity under fibrewise real analyticity in dimensions at least four and studied the associated twisted Codazzi quartics. The smooth angular quartic result used here is Proposition~\ref{prop:angular-rigidity}. Mao's recent preprint~\cite{MaoSurface} constructs a regular non-Berwald Landsberg surface on an open planar set, a two-dimensional noncompact example independent of Theorem~\ref{thm:main}.

In a fibrewise averaging argument, the differential of nonlinear parallel translation depends on the fibre variable. Indeed, if $P$ denotes nonlinear parallel translation, then
\[
 g_y(u,v)=g_{P(y)}\bigl((\dd P)_yu,(\dd P)_yv\bigr),
\]
so $(\dd P)_y$ cannot be replaced by a fixed linear map when the identity is integrated over an indicatrix. This issue is discussed in Matveev~\cite{Matveev} and Szab\'o~\cite{SzaboCorrection}. The argument below instead establishes linearity of the curvature fields in the original tangent-space coordinates before compact rigidity is invoked.

Fix $x\in M$, put $V=T_xM$, and let $(\Sigma,h)$ be the indicatrix in $V$, equipped with its intrinsic Riemannian metric. Consider the compact connected group
\[
K=\Isom(\Sigma,h)^0.
\]
The nonlinear curvature gives a family of Killing fields $R(u,v)$ on $\Sigma$, indexed by fixed vectors $u,v\in V$. The vertical Bianchi identity relates their derivatives in the linear fibre variable. The proof divides into two cases, according to whether $K$ acts transitively on $\Sigma$.

If $K$ is nontransitive, every curvature field is tangent to the orbits of $K$, and angular normal directions constrain the radial contraction of the curvature. The vertical Bianchi identity, the second fundamental form of an orbit, and its invariance under isotropy imply that $R(y,v)$ vanishes whenever $v$ is angular and normal to the orbit through $y$; see Theorem~\ref{thm:normal-vanishing}. The argument uses the isometry orbits themselves, particularly their invariant second fundamental forms. The resulting restriction on the fixed coefficient two-forms, together with the vertical Bianchi identity, implies that every $R(u,v)$ is linear in the original fibre variable.

If $K$ is transitive, there are no angular normal directions, but transitivity propagates any intrinsic curvature identity from one point to the entire indicatrix. If the Berwald tensor is nonzero on a fibre, Proposition~\ref{prop:angular-rigidity} yields a point of curvature isotropy on its indicatrix. Transitivity and Lemma~\ref{lem:curvature-normalization} then give constant sectional curvature one. Nonlinear parallel translation carries this property to every indicatrix in the connected component, and Theorem~\ref{prop:round-linearity} proves linearity of the curvature fields on each fibre.

Proposition~\ref{prop:angular-rigidity} is proved in Appendix~\ref{app:angular-proof}. The rank dichotomy follows the general framework of Li--Liu~\cite{LiLiu}, but the pointwise argument is different: Lemmas~\ref{alg:quartic-compatible} and~\ref{alg:jordan-closure} establish curvature compatibility and Jordan closure directly from the curvature symmetries and the Bianchi identity. In the transitive case, one isotropic point suffices, since transitivity propagates the intrinsic curvature identity.

On each connected component of $M$, the fibrewise results imply $R$-quadraticity, and Crampin's theorem completes the proof. Compactness of the base is needed only for this last implication. The same framework also recovers closed three-dimensional Landsberg rigidity. It does not, however, address the Laugwitz conjecture: Theorem~\ref{prop:round-linearity} proves linearity of the curvature fields, not vanishing of the Cartan cubic. Further dimension and regularity comments are collected in Appendix~\ref{app:variants}.

\begin{remark}
The conclusion of Theorem~\ref{thm:main} does not imply that $F$ is Riemannian. On $\R^n$, consider the translation-invariant Randers norm
\[
F(y)=\sqrt{\sum_i(y^i)^2}+b y^1,\qquad 0<|b|<1.
\]
This norm is regular and strongly convex, and it induces a Finsler metric on $\R^n/\mathbb Z^n$. Its energy is independent of $x$, so $G^i=0$ and the metric is Berwald. Since $b\neq0$, it is nonreversible and non-Riemannian.
\end{remark}

\section*{Acknowledgement}

The authors gratefully acknowledge the support of the Fundamental and Interdisciplinary Disciplines Breakthrough Plan of the Ministry of Education of China (JYB2025XDXM112).

\section{Finsler metrics and Hessian indicatrices}\label{sec:preliminaries}

\begin{definition}[Positive homogeneity]\label{def:homogeneity}
A function or vector-valued map $T$ on a slit vector space $V^\circ=V\setminus\{0\}$ is \emph{positively homogeneous of degree $r$} if
\[
T(ty)=t^rT(y)\qquad(t>0,\ y\in V^\circ).
\]
For tensor components on a tangent fibre, this identity is interpreted in fixed linear coordinates, with all argument vectors held fixed.
\end{definition}

\begin{definition}[Minkowski norm and regular Finsler metric]\label{def:Minkowski}
A \emph{smooth strongly convex Minkowski norm} on a finite-dimensional real vector space $V$ is a continuous function $F:V\to[0,\infty)$ that vanishes at zero, is positive and smooth on $V^\circ$, is positively homogeneous of degree one, and has positive-definite Hessian $D^2(F^2/2)$ on $V^\circ$. Reversibility $F(-y)=F(y)$ is not part of this definition.

A \emph{regular smooth strongly convex Finsler metric} on $M$ is a continuous function $F:TM\to[0,\infty)$, smooth on $T^\circ M$, whose restriction to every tangent space is a Minkowski norm in this sense. Regularity concerns all nonzero tangent directions.
\end{definition}

Let $\pi:T^\circ M\to M$ be the bundle projection. Its pullback tangent bundle is
\[
\pi^*TM=\{(x,y,v):y,v\in T_xM,\ y\ne0\}\longrightarrow T^\circ M,
\qquad (x,y,v)\longmapsto(x,y).
\]
Thus the fibre of $\pi^*TM$ over $(x,y)$ is $T_xM$. A tensor on this pullback bundle assigns a multilinear form on $T_xM$ to each nonzero reference vector $y\in T_xM$. Fibre differentiation uses the canonical identification $T_y(T_xM)=T_xM$.

In local coordinates $(x^i,y^i)$ on $TM$, put
\begin{equation}\label{eq:energy}
 E=\tfrac12 F^2,\qquad g_{ij}=\frac{\partial^2E}{\partial y^i\partial y^j},\qquad
 C_{ijk}=\tfrac12\frac{\partial g_{ij}}{\partial y^k}.
\end{equation}
The function $E$ is the \emph{energy}, $g$ is the \emph{fundamental tensor}, and $C$ is the \emph{Cartan cubic}, also called the Cartan tensor. Thus $g$ and $C$ are sections of $\Sym^2(\pi^*T^*M)$ and $\Sym^3(\pi^*T^*M)$, respectively. For every $y\in T_xM\setminus0$, the form $g_y$ is an inner product on $T_xM$. The condition $C=0$ is equivalent to the metric being Riemannian.

Euler's identity and its first two derivatives give
\begin{equation}\label{eq:euler}
 E_{y^i}y^i=2E,\qquad g_{ij}y^j=E_{y^i},\qquad
 g_y(y,y)=2E(y),\qquad C_y(y,\cdot,\cdot)=0.
\end{equation}
Moreover, $g_{ty}=g_y$ and $C_{ty}=t^{-1}C_y$ for $t>0$ in a fixed linear fibre coordinate system. The last equation in~\eqref{eq:euler} is called radial nullity.

On a fixed fibre $V$, let $D$ be its standard flat connection. We write $C_u$ for the endomorphism determined by
\[
g_y(C_uv,w)=C_y(u,v,w).
\]
Each $C_u$ is self-adjoint, and $C_uv=C_vu$. When the reference vector needs to be displayed, we also use $C^\sharp_y(u)v=C_uv$ and $C^\sharp_y(u,v)=C_uv$. For endomorphisms $A,B$, the commutator is $[A,B]=AB-BA$.

\begin{definition}[Hessian metric]\label{def:Hessian}
A Riemannian metric $g$ on a manifold with a flat torsion-free connection $D$ is a Hessian metric if locally $g=D\dd\varphi$ for a smooth function $\varphi$. In $D$-affine coordinates, this means $g_{ij}=\partial_i\partial_j\varphi$. For a Minkowski norm, $g=D^2E$ is a Hessian metric on $V^\circ$, with global potential $E=F^2/2$.
\end{definition}

\begin{lemma}\label{lem:hessianconnection}
The Levi--Civita connection of the Riemannian metric $g=D^2E$ on $V^\circ=V\setminus0$ is
\[
\nabla^g=D+C^\sharp.
\]
With the curvature convention $R(X,Y)Z=\nabla_X\nabla_YZ-\nabla_Y\nabla_XZ-\nabla_{[X,Y]}Z$, its curvature is
\begin{equation}\label{eq:hessian-curvature}
R^g(u,v)=-[C^\sharp(u),C^\sharp(v)].
\end{equation}
\end{lemma}
\begin{proof}
For constant vector fields the Koszul formula gives
\[
2g(\nabla^g_uv,w)=D_ug(v,w)+D_vg(u,w)-D_wg(u,v)=2C(u,v,w).
\]
Thus the connection matrix is $\Gamma=\frac12g^{-1}\dd g$. Set $A=g^{-1}\dd g$. Differentiating the inverse matrix gives $\dd A(u,v)=-[A(u),A(v)]$, so
\[
\dd\Gamma(u,v)=-2[\Gamma(u),\Gamma(v)].
\]
The curvature matrix $\dd\Gamma+[\Gamma,\Gamma]$ is therefore $-[\Gamma,\Gamma]$, which is~\eqref{eq:hessian-curvature}. All derivatives here are in the fixed linear coordinates of the same fibre.
\end{proof}

\begin{definition}
The \emph{indicatrix} of a Minkowski norm and its \emph{indicatrix metric} are
\[
\Sigma=\{y\in V:F(y)=1\},\qquad h_y=g_y|_{T_y\Sigma}.
\]
The metric $h$ is also called the angular metric. An \emph{angular vector} at $y\in\Sigma$ is an element of $T_y\Sigma=\ker\dd F_y=\ker\dd E_y$. The angular tensor $g-\dd F\otimes\dd F$ on $V^\circ$ restricts to $h$ on $T\Sigma$.
\end{definition}

Choose any auxiliary Euclidean norm $|\cdot|$ on $V$. Positivity and continuity on its unit sphere give constants $0<c\leq C_0<\infty$ such that $c|y|\leq F(y)\leq C_0|y|$. The radial map $u\mapsto u/F(u)$ is a diffeomorphism from the Euclidean unit sphere onto $\Sigma$. Thus $\Sigma$ is compact, connected, and diffeomorphic to $S^{n-1}$. It is simply connected when $n\geq3$.

Write $y=r\theta$, where $r=F(y)>0$ and $\theta\in\Sigma$. For $u\in T_\theta\Sigma$, Euler's identity gives $g_\theta(\theta,u)=0$ and $g_\theta(\theta,\theta)=1$. Consequently
\begin{equation}\label{eq:cone}
g=\dd r^2+r^2h.
\end{equation}
The cone connection at $r=1$ satisfies
\[
\nabla^g_uv=\nabla^h_uv-h(u,v)y
\]
for angular tangent fields. Combining this with Lemma~\ref{lem:hessianconnection} gives
\begin{equation}\label{eq:angularD}
D_uX=\nabla^h_uX-C^\sharp(u,X)-h(u,X)y
\end{equation}
for a field tangent to the indicatrix and extended as a positively homogeneous field of degree one. Its radial derivative is $D_yX=X$.

Define $R_0(u,v)w=h(v,w)u-h(u,w)v$. The cone curvature and~\eqref{eq:hessian-curvature} give the angular Gauss equation
\begin{equation}\label{eq:Gauss}
R^h(u,v)=R_0(u,v)-[C^\sharp(u),C^\sharp(v)].
\end{equation}
\begin{definition}[Angular sectional curvature]\label{def:angular-curvature}
For linearly independent $u,v\in T_y\Sigma$, the sectional curvature of their span is
\[
\sec_h(u\wedge v)=
\frac{h(R^h(u,v)v,u)}{h(u,u)h(v,v)-h(u,v)^2}.
\]
The angular curvature is \emph{isotropic at $y$} if $R^h_y=\kappa(y)R_{0,y}$ for a scalar $\kappa(y)$. It has \emph{constant sectional curvature $\kappa$} if this equality holds with the same scalar at every point. An indicatrix is \emph{round with sectional curvature one} if it is isometric to the Euclidean unit sphere with its induced metric. In the compact simply connected setting considered here, constant sectional curvature one is equivalent to this condition.
\end{definition}

The restriction of $C$ to $T\Sigma$ determines the full cubic by radial nullity. It is a \emph{Codazzi cubic}: the covariant four-tensor $\nabla^hC$ is totally symmetric. To verify this, use affine fibre coordinates and write $E_{ijkl}$ for the fourth derivative of $E$. The connection in Lemma~\ref{lem:hessianconnection} gives
\[
(\nabla^g_\ell C)_{ijk}
=\frac12E_{ijkl}
-g^{pq}\bigl(C_{l ip}C_{jkq}+C_{l jp}C_{ikq}+C_{l kp}C_{ijq}\bigr).
\]
The right-hand side is symmetric in $i,j,k,\ell$. Its restriction to angular arguments is $\nabla^hC$, since the radial terms in the cone connection have zero insertion in $C$. A symmetric cubic is called \emph{commuting} when its metric-dual endomorphisms satisfy $[C_u,C_v]=0$ for all $u,v$ at every point.

The tensor $R^g$ in~\eqref{eq:hessian-curvature} is the Riemannian curvature of the Hessian metric on one punctured tangent space. It is distinct from the curvature of the horizontal distribution over the base, represented by the fields $R(u,v)$ below.

The canonical spray coefficients are
\begin{equation}\label{eq:spray}
G^m=\frac12g^{mi}\bigl(E_{x^ky^i}y^k-E_{x^i}\bigr).
\end{equation}
The associated \emph{canonical spray} is the vector field
\[
S=y^i\partial_{x^i}-2G^m\partial_{y^m}
\]
on $T^\circ M$. Its integral curves project to geodesics satisfying $\ddot x^m+2G^m(x,\dot x)=0$. The coefficients $G^m$ are positively homogeneous of degree two. Put
\[
N^m_i=\partial_{y^i}G^m,\qquad \delta_i=\partial_{x^i}-N^m_i\partial_{y^m},\qquad
\Gamma^m_{ij}=\partial_{y^j}N^m_i.
\]
The fields $\delta_i$ span the horizontal distribution of the canonical nonlinear connection; the vertical distribution is $\ker\dd\pi$, spanned by $\partial_{y^i}$. The symmetry $\Gamma^m_{ij}=\Gamma^m_{ji}$ follows from the fibre Hessian of $G^m$.

The \emph{Berwald connection} $\nabla^B$ on $\pi^*TM$ is determined in the induced frame $e_j=\partial_{x^j}$ by
\[
\nabla^B_{\delta_i}e_j=\Gamma^m_{ij}e_m,\qquad
\nabla^B_{\partial_{y^i}}e_j=0.
\]
A horizontal derivative of a covariant tensor differentiates its components with $\delta_i$ and subtracts the corresponding connection term in each covariant slot. We denote it by a vertical bar.

The Berwald curvature is
\begin{equation}\label{eq:berwald}
B^m{}_{ijk}=\partial_{y^i}\partial_{y^j}\partial_{y^k}G^m.
\end{equation}
It is symmetric in its three lower indices and has zero contraction with $y$ in any of them. The metric is Berwald precisely when this tensor is zero. Equivalently, the spray is quadratic in $y$, with coefficients depending only on $x$.

The \emph{Landsberg tensor} is $L_{ijk}=C_{ijk|l}y^l$. A Finsler metric is called \emph{Landsberg} if $L=0$. In terms of the Berwald curvature this condition is equivalent to
\begin{equation}\label{eq:landsberg}
E_{y^m}B^m{}_{ijk}=0.
\end{equation}
The canonical nonlinear connection is conservative:
\begin{equation}\label{eq:conservative}
\delta_iE=0.
\end{equation}
Differentiating twice in the fibre gives
\begin{equation}\label{eq:horizontalmetric}
Q_{ijk}:=\delta_i g_{jk}-\Gamma^m_{ij}g_{mk}-\Gamma^m_{ik}g_{jm}
=E_{y^m}B^m{}_{ijk}.
\end{equation}
Differentiating~\eqref{eq:conservative} once in the fibre gives $\delta_i E_{y^j}=(\partial_{y^j}N^m_i)E_{y^m}$. Differentiating this identity in $y^k$ gives the two connection terms and the right side of~\eqref{eq:horizontalmetric}. Consequently, the Landsberg condition is equivalent to horizontal metricity of the Berwald connection.

\section{Landsberg transport and curvature fields}\label{sec:curvature}

We use the nonlinear-curvature convention of Mao--Zhou~\cite{MaoZhou}:
\begin{equation}\label{eq:rawR}
R^m{}_{ij}=\delta_jN^m_i-\delta_iN^m_j.
\end{equation}
Thus $[\delta_i,\delta_j]=R^m{}_{ij}\partial_{y^m}$. At a fixed base point $x$ and for vectors $u,v\in V=T_xM$, define the \emph{curvature field}
\[
R_x(u,v)(y)=u^iv^jR^m{}_{ij}(x,y)\partial_{y^m}.
\]
When $x$ is fixed, we abbreviate $R_x(u,v)$ to $R(u,v)$. This dependence on $u,v$ is alternating and bilinear. The arguments $u,v$ are fixed elements of $V$, while $y$ is the variable of the vector field on $V^\circ$.

\begin{definition}[$R$-quadraticity]\label{def:R-quadratic}
The \emph{Jacobi endomorphism}, with the curvature sign used in this paper, is the family of linear maps
\[
\rho_{x,y}:T_xM\longrightarrow T_xM,\qquad
\rho_{x,y}(v)=R(y,v)(y).
\]
It is positively homogeneous of degree two in the reference vector $y$. A Finsler metric is \emph{$R$-quadratic} if the components of $\rho_{x,y}$ are homogeneous polynomials of degree two in the linear fibre coordinates, with coefficients depending smoothly on $x$. The equivalence with linearity of $R^m{}_{ij}(x,y)$ in $y$ is established in Section~\ref{sec:compact}.
\end{definition}

\begin{definition}[Killing fields and isometry orbits]\label{def:Killing}
On a Riemannian manifold $(N,h)$, a \emph{Killing field} is a smooth vector field $X$ satisfying $\Lie_Xh=0$, or equivalently
\[
h(\nabla^h_UX,V)+h(U,\nabla^h_VX)=0
\]
for all tangent vectors $U,V$. Its local flow consists of isometries. We denote the Lie algebra of global Killing fields, with the usual vector-field bracket, by $\Kill(N,h)$.

For a group $K$ acting on $N$, the \emph{orbit} of $p$ is $K\cdot p=\{\phi(p):\phi\in K\}$. The action is \emph{transitive} if this orbit is $N$. Otherwise it is \emph{nontransitive}. In this paper $K=\Isom(\Sigma,h)^0$, the connected component of the identity in the intrinsic isometry group of the indicatrix.
\end{definition}

\begin{lemma}\label{lem:curvature-identities}
The fields $R(u,v)$ are smooth on $V^\circ$, positively one-homogeneous, tangent to every level of $E$, and Killing for $g$. They satisfy
\begin{equation}\label{eq:VBianchi}
D_uR(v,w)+D_vR(w,u)+D_wR(u,v)=0.
\end{equation}
\end{lemma}
\begin{proof}
The degrees of $N$ and $R$ follow directly from the degree of $G$. Taking two horizontal derivatives of~\eqref{eq:conservative} in either order gives $\dd E(R(u,v))=0$.

The expression $Q_{ijk}$ in~\eqref{eq:horizontalmetric} is the Lie derivative of the vertical metric under $\delta_i$, evaluated on $\partial_{y^j},\partial_{y^k}$. Indeed, $[\delta_i,\partial_{y^j}]=\Gamma^m_{ij}\partial_{y^m}$. The horizontal flows preserve the vertical bundle, and the Landsberg condition gives $\Lie_{\delta_i}g=0$ on that bundle. Consequently
\[
\Lie_{[\delta_i,\delta_j]}g=[\Lie_{\delta_i},\Lie_{\delta_j}]g=0.
\]
The bracket is the curvature field $R_{ij}$, so $\Lie_{R(u,v)}g=0$.

Fibre differentiation of~\eqref{eq:rawR} gives
\[
\partial_{y^k}R^m{}_{ij}
=\delta_j\Gamma^m_{ik}-\delta_i\Gamma^m_{jk}
-\Gamma^p_{jk}\Gamma^m_{ip}+\Gamma^p_{ik}\Gamma^m_{jp}.
\]
Since $\Gamma^m_{ij}=\Gamma^m_{ji}$, the cyclic sum in $i,j,k$ vanishes term by term. This is~\eqref{eq:VBianchi} for coordinate labels, and multilinearity gives the result for all fixed $u,v,w\in V$.
\end{proof}

When a label is written as $y$ at a specified point, it is held fixed during differentiation. Thus $(D_wR(y,v))(y)$ differentiates the field whose first label is that fixed vector, whereas differentiation of $z\mapsto R(z,v)(z)$ also differentiates the first label. All applications of the vertical Bianchi identity use fixed labels.

Let $\gamma:[0,1]\to M$ be a smooth curve, written $x(t)$ in coordinates. For $y_0\in T_{\gamma(0)}M\setminus0$, its \emph{nonlinear parallel lift} is the solution of
\[
\dot y^m+N^m_i(x(t),y)\dot x^i=0,\qquad y(0)=y_0.
\]
Conservativity gives $F(x(t),y(t))=F(\gamma(0),y_0)$. On the compact parameter interval, the corresponding indicatrices form a compact bundle, so the lift exists up to $t=1$. The endpoint map $P_\gamma(y_0)=y(1)$ is called \emph{nonlinear parallel translation}. Reversing $\gamma$ gives its inverse. Homogeneity makes $P_\gamma$ a diffeomorphism between the slit tangent fibres that commutes with positive dilations. For a piecewise smooth curve, these maps are composed over its smooth subintervals.

Differentiating the lift equation with respect to its initial vector shows that a variation $v$ satisfies
\[
\dot v^m+\Gamma^m_{ij}\dot x^i v^j=0.
\]
Equation~\eqref{eq:horizontalmetric} now gives
\[
\frac{\dd}{\dd t}g_{y(t)}(v(t),w(t))
=Q_{ijk}\dot x^iv^jw^k=0.
\]
Thus $P_\gamma$ restricts to an isometry $(\Sigma_{\gamma(0)},h_{\gamma(0)})\to(\Sigma_{\gamma(1)},h_{\gamma(1)})$. Push-forward identifies their Killing algebras, and conjugation identifies the identity components of their isometry groups. In particular, both $\dim\Kill(\Sigma_x,h_x)$ and transitivity of $\Isom(\Sigma_x,h_x)^0$ are constant on each connected component of the base.

Every Killing field on the compact Riemannian manifold $(\Sigma,h)$ is complete, and its flow belongs to $K$. Thus each $R(u,v)|_\Sigma$ is tangent to every $K$-orbit. If every orbit is a point, all curvature fields vanish. Otherwise there is a positive-dimensional orbit.

\begin{remark}
An intrinsic isometry preserves the Riemannian metric $h$ on $\Sigma$. A linear norm isometry is a linear automorphism $A$ of $V$ satisfying $F\circ A=F$; its restriction to $\Sigma$ is an intrinsic isometry. No extension of an arbitrary intrinsic isometry to a linear map of $V$ is assumed.
\end{remark}

\section{Isometry orbits and normal curvature fields}\label{sec:shape}

\begin{definition}[Intrinsic isometry action]\label{def:isometry-action}
The group $K=\Isom(\Sigma,h)^0$ acts on $\Sigma$ by evaluation,
$(k,z)\mapsto k(z)$. Its orbit through $y$ is
$K\cdot y=\{k(y):k\in K\}$, and its isotropy subgroup at $y$ is
$K_y=\{k\in K:k(y)=y\}$. The differential action
$k\mapsto (\dd k)_y$ of $K_y$ on $T_y\Sigma$ is its isotropy representation.
An orbit is proper if it is strictly contained in $\Sigma$.
\end{definition}

\begin{definition}[Fixed-label curvature family]\label{def:fixed-label-family}
A family of curvature fields with fixed arguments is an alternating bilinear map
$V\times V\to\mathfrak X(V^\circ)$, $(u,v)\mapsto R(u,v)$, where
$\mathfrak X(V^\circ)$ denotes the space of smooth vector fields on
$V^\circ$. The labels $u,v\in V$ are fixed vectors, independent of the variable of
the resulting field, and all fibre derivatives use the fixed-label convention
introduced in Section~\ref{sec:curvature}. Coordinate tensor indices are
written as $i,j,k,l,m,p,q$.
\end{definition}

Fix one fibre and let $Y=K\cdot y$ be an orbit in $\Sigma$. Compact group orbits are compact embedded submanifolds, and $Y$ is connected because $K$ is connected. At the reference point put
\[
T=T_yY,\qquad N=T^{\perp_h}\cap T_y\Sigma,\qquad V=\R y\oplus T\oplus N.
\]
The sum is orthogonal for $g_y$. In each fibre derivative the arguments $u,v$ of $R(u,v)$ are held fixed in $V$. Since $R(u,v)|_\Sigma$ is a global Killing field, its value $R(u,v)(y)$ belongs to $T$.

\begin{theorem}[Normal-label vanishing]\label{thm:normal-vanishing}
For $\nu\in N$, the fixed-label field $R(y,\nu)$ vanishes identically on $V^\circ$.
\end{theorem}
We first derive the identities for the second fundamental form needed in the proof.

For an angular field $X$, equation~\eqref{eq:angularD} implies
\begin{equation}\label{eq:radialD}
g_y(D_uX,y)=-g_y(u,X(y)),\qquad u\in T_y\Sigma.
\end{equation}
Apply the vertical Bianchi identity to $(y,u,v)$, with $u,v$ angular. Its radial component gives
\[
g_y\bigl(u,R(y,v)(y)\bigr)=g_y\bigl(v,R(y,u)(y)\bigr).
\]
The map $u\mapsto R(y,u)(y)$ is therefore self-adjoint on the angular tangent space. Its range is in $T$, so it annihilates $N$. Write
\begin{equation}\label{eq:Mdef}
M:T\longrightarrow T,\qquad Mu=R(y,u)(y).
\end{equation}
Thus $M$ is self-adjoint with respect to $h_y$.

For $\nu\in N$ set $S_\nu=R(y,\nu)|_\Sigma$. Its value at $y$ is zero. Its isometric flow fixes $y$ and preserves $Y$, so
\[
A_\nu=(\nabla^h S_\nu)_y
\]
is skew-adjoint and preserves $T$ and $N$.

\begin{lemma}\label{lem:A-normal}
For every $\nu\in N$, $A_\nu|_N=0$.
\end{lemma}
\begin{proof}
Apply the vertical Bianchi identity to $(y,\nu,\nu')$, where $\nu,\nu'\in N$, and take its $N$ component. The radial derivative term is tangent to $Y$. In the other terms the fields vanish at $y$, so the Cartan and radial terms in~\eqref{eq:angularD} vanish. Hence $A_\nu\nu'=A_{\nu'}\nu$. The tensor
\[
\mathfrak a(\nu,\nu',\nu'')=g_y(A_\nu\nu',\nu'')
\]
on $N$ is symmetric in its first two slots and skew in its last two. Explicitly,
\begin{align*}
\mathfrak a(\nu,\nu',\nu'')
&=\mathfrak a(\nu',\nu,\nu'')
=-\mathfrak a(\nu',\nu'',\nu)
=-\mathfrak a(\nu'',\nu',\nu)\\
&=\mathfrak a(\nu'',\nu,\nu')
=\mathfrak a(\nu,\nu'',\nu')
=-\mathfrak a(\nu,\nu',\nu'').
\end{align*}
Consequently $\mathfrak a=0$. Since $A_\nu$ preserves $N$, it follows that $A_\nu|_N=0$.
\end{proof}

Define $L_\nu:T\to T$ by $L_\nu u=R(\nu,u)(y)$. The radial component of the vertical Bianchi identity for $(\nu,u,v)$, with $u,v\in T$, gives
\[
g_y\bigl(u,R(\nu,v)(y)\bigr)=g_y\bigl(v,R(\nu,u)(y)\bigr).
\]
The term paired with $\nu$ is zero because all curvature values are tangent to $Y$. Thus $L_\nu$ is self-adjoint.

\begin{definition}[Second fundamental form and shape operator]\label{def:orbit-shape}
For tangent fields $X,Z$ along $Y$, the second fundamental form of
$Y\subset(\Sigma,h)$ is the normal-valued symmetric tensor
$\mathrm{II}(X,Z)=(\nabla^h_XZ)^N$. Here the normal bundle is the
$h$-orthogonal complement of $TY$ in $T\Sigma|_Y$. For a normal vector
$\nu\in N$, the associated shape operator $\mathrm{II}_\nu:T\to T$ is defined by
$h(\mathrm{II}_\nu u,v)=h(\nu,\mathrm{II}(u,v))$. This convention determines its sign.
\end{definition}

At the reference point define the symmetric Cartan contraction $C_\nu:T\to T$ by
\[
g_y(C_\nu u,v)=C_y(\nu,u,v).
\]
Since $R(y,u)$ is tangent to $Y$ everywhere on the orbit and has value $Mu$ at $y$,
\[
(\nabla^h_vR(y,u))^N=\mathrm{II}(v,Mu).
\]
The Killing equation then gives
\begin{equation}\label{eq:shape-killing}
(\nabla^h_\nu R(y,u))^T=-\mathrm{II}_\nu Mu.
\end{equation}
Here the projections are taken at the reference point.

The vertical Bianchi identity with arguments $(\nu,y,u)$ reads
\[
D_\nu R(y,u)+D_yR(u,\nu)+D_uR(\nu,y)=0.
\]
Its tangential parts are respectively $-(\mathrm{II}_\nu+C_\nu)Mu$, $-L_\nu u$, and $-A_\nu u$. The first follows from~\eqref{eq:angularD} and~\eqref{eq:shape-killing}, the second from homogeneity, and the third from $R(\nu,y)|_\Sigma=-S_\nu$ and $S_\nu(y)=0$. Hence
\begin{equation}\label{eq:shape-main}
A_\nu|_T+L_\nu=-(\mathrm{II}_\nu+C_\nu)M.
\end{equation}
Taking skew parts gives
\begin{equation}\label{eq:shape-skew}
A_\nu|_T=-\tfrac12[\mathrm{II}_\nu+C_\nu,M].
\end{equation}

Next apply the identity to $(u,y,v)$, take its normal component, and pair with $\nu$. The radial derivative term is tangent to $Y$. The remaining terms give
\[
g_y((\mathrm{II}_\nu-C_\nu)u,Mv)=g_y((\mathrm{II}_\nu-C_\nu)v,Mu).
\]
Since both endomorphisms are symmetric,
\begin{equation}\label{eq:shape-othercomm}
[\mathrm{II}_\nu-C_\nu,M]=0.
\end{equation}
Combining the two commutator identities yields
\begin{equation}\label{eq:Acomm}
A_\nu|_T=-[\mathrm{II}_\nu,M].
\end{equation}

\begin{proof}[Proof of Theorem~\ref{thm:normal-vanishing}]
Let $A=A_\nu|_T$. The flow of $S_\nu$ fixes $y$ and preserves the orbit and its second fundamental form. Its derivative on $N$ is the identity, because the generator there is zero by Lemma~\ref{lem:A-normal}. Differentiating the invariance of the second fundamental form gives
\[
\mathrm{II}(Au,v)+\mathrm{II}(u,Av)=0.
\]
Pairing with $\nu$ and using skewness of $A$ gives
\begin{equation}\label{eq:shape-isotropy}
[A,\mathrm{II}_\nu]=0.
\end{equation}
Thus cyclicity of the trace and~\eqref{eq:Acomm} show
\begin{align*}
\norm{A}_{\mathrm{HS}}^2
&=-\tr(A^2)=\tr\bigl(A[\mathrm{II}_\nu,M]\bigr)\\
&=\tr\bigl([A,\mathrm{II}_\nu]M\bigr)=0.
\end{align*}
Here $\norm{A}_{\mathrm{HS}}^2=\tr(A^*A)$ is the squared Hilbert--Schmidt
norm, and $A^*=-A$ because $A$ is skew-adjoint.
The full angular derivative $A_\nu$ vanishes.

A global Killing field on a connected Riemannian manifold is determined by its value and covariant derivative at one point. Indeed, its restriction to a geodesic is a Jacobi field. Vanishing initial data therefore imply vanishing along that geodesic, and connectedness extends the conclusion to the manifold. Since $S_\nu(y)=0$ and $(\nabla^hS_\nu)_y=0$, we obtain $S_\nu=0$ on $\Sigma$. Positive homogeneity gives $R(y,\nu)=0$ on $V^\circ$.
\end{proof}

\begin{remark}
The invariance of $\mathrm{II}$ used in \eqref{eq:shape-isotropy} follows from its definition as the second fundamental form of an isometry orbit. It is this invariance that permits the trace identity to imply vanishing of the commutator.
\end{remark}

\section{Curvature fields in the nontransitive case}\label{sec:proper}

Assume that $K$ is nontransitive. If every curvature field vanishes, linearity is immediate. Otherwise choose a positive-dimensional orbit $Y$. It has dimension strictly less than $n-1$: a full-dimensional compact orbit would be open and closed in the connected indicatrix and hence equal to it.

\begin{theorem}\label{thm:proper-linearity}
If $K$ is nontransitive, each curvature field $R(u,v)$ is linear in the original fibre coordinates.
\end{theorem}

Consider $\mathcal R:\Lambda^2V\to\mathfrak X(V^\circ)$, $u\wedge v\mapsto R(u,v)$, and set
\begin{equation}\label{eq:Wdef}
\mathcal W=(\ker\mathcal R)^\perp\subset\Lambda^2V^*.
\end{equation}
For a subspace $L\subset\Lambda^2V$, its annihilator is
$L^\perp=\{\omega\in\Lambda^2V^*: \omega(\ell)=0\text{ for every }\ell\in L\}$;
the superscript in \eqref{eq:Wdef} has this meaning rather than denoting an
inner-product orthogonal complement. Equivalently, $\mathcal W$ is the span
of all two-forms obtained by evaluating the coordinate components of the
vector-valued curvature $R$ at arbitrary points. Indeed a bivector is
annihilated by all those evaluations precisely when its associated curvature
field is zero; finite-dimensional duality gives the equality of the two definitions.

For $\omega\in\mathcal W$ put
\begin{equation}\label{eq:Xomega}
X_\omega(y)=g_y^{-1}\bigl(\omega(y,\cdot)\bigr).
\end{equation}
Alternation makes this field angular. If $y\in Y$ and $\nu$ is angular-normal, Theorem~\ref{thm:normal-vanishing} gives $y\wedge\nu\in\ker\mathcal R$, hence $\omega(y,\nu)=0$. Consequently $X_\omega$ is tangent to $Y$. It need not be Killing. At its zeros, its differential has the skew-adjoint form used below.

\begin{lemma}\label{lem:good-symbol}
There is $\xi\in V^*\setminus0$ with
\begin{equation}\label{eq:good-symbol}
\mathcal W\cap(\xi\wedge V^*)=\{0\}.
\end{equation}
\end{lemma}
\begin{proof}
Let $P=\operatorname{span}(Y)$ and $d=\dim Y$. Since $\R y+T_yY\subset P$, $\dim P\geq d+1$.

\emph{First case: $\dim P=d+1$.}
The inclusion $Y\subset\Sigma\cap P$ is between submanifolds of the same dimension. It is open by the inverse function theorem and closed by compactness. The radial section $\Sigma\cap P$ is a connected sphere of dimension $d\geq1$, so
\[
Y=\Sigma\cap P.
\]
The subspace $P$ is proper because the orbit is proper. Choose $0\ne\delta\in\Ann(P)$. Its $g_y$-dual is angular-normal and nonzero along $Y$. The norm of its normal projection has a positive minimum on $Y$. Choose $\xi$ close to $\delta$, with $\xi|_P\ne0$, so that $g_y^{-1}\xi$ still has nonzero normal projection everywhere on $Y$.

Suppose $\xi\wedge\eta\in\mathcal W$. At $y\in Y\cap\ker\xi$,
\[
X_{\xi\wedge\eta}(y)=-\eta(y)g_y^{-1}\xi.
\]
Tangency forces $\eta(y)=0$. This section spans $P\cap\ker\xi$: each nonzero vector in that subspace has a positive radial rescaling in the section. Thus $\eta|_P=\lambda\xi|_P$ for some $\lambda\in\R$.

Set $\eta'=\eta-\lambda\xi\in\Ann(P)$. At a point $y\in Y$ with $\xi(y)\ne0$,
\[
X_{\xi\wedge\eta'}(y)=\xi(y)g_y^{-1}\eta'
\]
is tangent to $Y$. But $g_y^{-1}\eta'$ is angular-normal there, so it must vanish. Thus $\eta'=0$ as a covector, and $\xi\wedge\eta=0$.

\emph{Second case: $\dim P>d+1$.}
Choose $y_0\in Y$ and a covector $\xi$ annihilating $\R y_0+T_{y_0}Y$ but not $P$. Suppose $\omega=\xi\wedge\eta\in\mathcal W$ is nonzero. At $y_0$,
\[
X_\omega(y_0)=-\eta(y_0)g_{y_0}^{-1}\xi
\]
is both normal and tangent, hence zero. Therefore $\eta(y_0)=0$. At any zero of $X_\omega$, the derivative of $g^{-1}$ contributes zero, and
\begin{equation}\label{eq:Xzero}
(DX_\omega)_yv=g_y^{-1}\bigl(\xi(v)\eta-\eta(v)\xi\bigr).
\end{equation}
It is a skew-adjoint rank-two map with image
\[
L_y=\operatorname{span}(g_y^{-1}\xi,g_y^{-1}\eta).
\]
At $y_0$ this plane is angular. The tangent space $T_{y_0}Y$ is invariant under the derivative because the field is tangent to $Y$ and vanishes at that point. The plane $L_{y_0}$ contains the nonzero normal vector $g^{-1}\xi$.

A rank-two skew-adjoint map acts as a nonzero rotation on its image plane and as zero on its orthogonal complement. Any invariant real subspace either contains this plane or is orthogonal to it. To see this, its square is a nonzero scalar times the projection onto the plane; applying that square to an invariant subspace gives its plane component, and a nonzero rotation has no invariant real line. The first possibility is excluded by the normal vector. Hence
\[
T_{y_0}Y\subset L_{y_0}^\perp,\qquad (DX_\omega)_{y_0}|_{T_{y_0}Y}=0.
\]

Let $\pi:V^*\to V^*/\operatorname{span}(\xi,\eta)$ and define
\[
f=\pi\circ\dd E:V^\circ\longrightarrow V^*/\operatorname{span}(\xi,\eta).
\]
Its derivative has kernel $L_y$ because $D(\dd E)_y=g_y$. Also
\[
Df_y(X_\omega(y))=\pi\bigl(\omega(y,\cdot)\bigr)=0.
\]
At $y_0$, $Df|_{TY}$ is injective. The immersion theorem therefore makes $f|_Y$ injective on a sufficiently small neighborhood. Along a short integral curve of $X_\omega$ in that neighborhood, $f$ is constant, so the curve is stationary. The same applies to every starting point in that neighborhood. Thus the field vanishes locally on $Y$.

Let $Z$ be the set of points with such a stationary neighborhood. It is nonempty and open. At a limit point of $Z$, continuity gives $X_\omega=0$ and $(DX_\omega)|_{TY}=0$. By~\eqref{eq:Xzero}, $TY\subset L_y^\perp$, so $Df|_{TY}$ is injective again. The same argument gives local stationarity. Thus $Z$ is closed. Connectedness implies $Z=Y$.

Independence of $\xi,\eta$ now implies $\xi(y)=\eta(y)=0$ for every $y\in Y$, so $\xi|_P=0$, a contradiction. Both cases prove the lemma.
\end{proof}

Fix linear coordinates on $V\simeq\R^n$. For each fixed component index $m$, let
\[
u(y)=\tfrac12R^m{}_{ij}(y)\,\dd y^i\wedge\dd y^j.
\]
This is a smooth map $V^\circ\to\mathcal W$, regarded as a two-form on $V^\circ$. The vertical Bianchi identity is precisely $\dd u=0$, where the exterior derivative is taken in linear coordinates on $V$.

The exterior multiplication map
$L_\xi:\mathcal W\to\Lambda^3V^*$, $L_\xi\omega=\xi\wedge\omega$, is the
symbol of this constant-coefficient first-order equation, up to the
nonzero scalar determined by the Fourier transform convention.
Two-form-valued distributions and their Fourier transforms are understood
componentwise in a fixed basis; a distribution is tempered if its
coefficients are continuous linear functionals on the Schwartz space of
smooth rapidly decreasing functions.

\begin{lemma}\label{lem:fourier-linearity}
Let $u:\R^n\setminus0\to\mathcal W\subset\Lambda^2(\R^n)^*$ be smooth and positively homogeneous of degree one. If $\dd u=0$ and~\eqref{eq:good-symbol} holds, then $u$ is linear.
\end{lemma}
\begin{proof}
Homogeneity and compactness of the Euclidean sphere give
\[
|D^\alpha u(y)|\leq C_\alpha |y|^{1-|\alpha|}.
\]
Thus $u$ is locally integrable across zero and has polynomial growth at infinity, defining a tempered distribution. The equation $\dd u=0$ holds across zero as well. Indeed, integrate by parts against a compactly supported smooth test form outside a ball of radius $\varepsilon$. The boundary integrand is $O(\varepsilon)$, while its area is $O(\varepsilon^{n-1})$. The error is $O(\varepsilon^n)$ and tends to zero.

Taking Fourier transforms gives
\begin{equation}\label{eq:Fourierclosed}
\zeta\wedge\widehat u(\zeta)=0
\end{equation}
as an identity of distributions; the nonzero scalar factor determined by the Fourier convention is omitted.

For $\xi\ne0$, the kernel of exterior multiplication by $\xi$ on two-forms is $\xi\wedge V^*$. Indeed, split $V^*=\R\xi\oplus U$ and decompose a two-form as $\xi\wedge\eta+\omega_U$. Exterior multiplication by $\xi$ is injective on $\Lambda^2U$. Therefore the symbol
\[
L_\xi:\mathcal W\longrightarrow\Lambda^3V^*,\qquad \omega\longmapsto\xi\wedge\omega
\]
is injective.

If $\mathcal W=0$, then $u=0$. Otherwise choose bases and a maximal minor of $L_\zeta$ that is nonzero at this $\xi$. Its determinant $p(\zeta)$ is a nonzero polynomial, because the matrix entries are linear in $\zeta$. On the open dense set $\{p\ne0\}$ the inverse minor gives a smooth left inverse, so~\eqref{eq:Fourierclosed} forces $\widehat u=0$ there as a distribution.

It remains to prove that $\widehat u$ is smooth away from zero. Choose a compactly supported smooth cutoff $\chi$ equal to one near zero and decompose
\[
u=\chi u+u_\infty .
\]
The compactly supported part is integrable with integrable moments of all orders, so its Fourier transform is smooth. For every multi-index $\alpha$,
\[
\partial_{y^j}^m(y^\alpha u_\infty)\in L^1(\R^n)
\quad\hbox{when }m>n+1+|\alpha|.
\]
There is no singularity near zero, and the derivative has decay $O(|y|^{|\alpha|+1-m})$ at infinity. On any neighborhood where $\zeta_j\ne0$, the Fourier derivative identity represents the distributional derivative $\partial_\zeta^\alpha\widehat u_\infty$ by a continuous Fourier transform of this integrable function divided by $(i\zeta_j)^m$, up to a constant. All distributional derivatives thus have continuous representatives locally. The usual weak-derivative regularity criterion yields smoothness. Every nonzero $\zeta$ has a nonzero coordinate, so $\widehat u$ is smooth off zero.

It vanishes on the dense open set $\{p\ne0\}$, hence by continuity everywhere off zero. A distribution supported at one point is a finite sum of derivatives of a delta distribution~\cite{Hormander}; Fourier inversion makes $u$ a polynomial. Finally the polynomial identity $u(ty)=tu(y)$, $t>0$, eliminates all degrees except one. Thus $u$ is linear.
\end{proof}

\begin{proof}[Proof of Theorem~\ref{thm:proper-linearity}]
If all curvature fields vanish, the assertion is immediate. Otherwise $K$ has a positive-dimensional proper orbit. Lemma~\ref{lem:good-symbol} gives a covector satisfying~\eqref{eq:good-symbol}. For each component index, the nonlinear curvature defines a map $u:V^\circ\to\mathcal W$ satisfying the hypotheses of Lemma~\ref{lem:fourier-linearity}. Its coefficients are therefore linear in $y$. This proves the assertion for every curvature field.
\end{proof}

\begin{remark}
Smoothness on the entire slit space gives the uniform derivative bounds on the Euclidean sphere used in the Fourier argument, and the global norm hypotheses give the compact radial sections used in the incidence argument. A conic or almost regular metric can fail these requirements. No antipodal symmetry is used: the argument requires only positive homogeneity and positive radial rescalings.
\end{remark}

\section{The transitive case and angular quartic rigidity}\label{sec:homogeneous}

Assume that $K=\Isom(\Sigma,h)^0$ acts transitively. The indicatrix isometries induced by nonlinear parallel translation imply that the same condition holds on every fibre over the connected component of the base under consideration. The action preserves the indicatrix metric $h$.

Let $B^m{}_{ijk}$ be the Berwald curvature and lower its vector index with the fundamental tensor:
\begin{equation}\label{eq:angularBerwaldDefinition}
B^m{}_{ijk}=\partial_{y^i}\partial_{y^j}\partial_{y^k}G^m,
\qquad B_{l ijk}=g_{l m}B^m{}_{ijk}.
\end{equation}
These tensors are defined on $\pi^*TM$. The covariant tensor $B_{l ijk}$ and its restriction to $T\Sigma_x$ will both be denoted by $B$; the arguments specify which restriction is intended.

\begin{lemma}\label{lem:angularnonzero}
For a Landsberg metric, $B_{l ijk}$ is completely symmetric, has degree $-1$ in the reference vector, and vanishes when the reference vector is inserted in any slot. A fibre has nonzero Berwald curvature if and only if its angular quartic is not identically zero.
\end{lemma}
\begin{proof}
The Landsberg equation is $E_mB^m{}_{ijk}=0$, where $E_m=\partial_{y^m}E$. Differentiating this scalar identity in the fibre direction $\partial_{y^l}$ gives
\[
g_{ml}B^m{}_{ijk}+E_m\partial_{y^l}\partial_{y^i}\partial_{y^j}\partial_{y^k}G^m=0.
\]
Consequently,
\begin{equation}\label{eq:symmetricB}
B_{l ijk}=-E_m\partial_{y^l}\partial_{y^i}\partial_{y^j}\partial_{y^k}G^m.
\end{equation}
The four fibre differentiations commute. This proves symmetry in all four slots. Since $G$ has degree two, its third derivatives have degree $-1$, whereas $g$ has degree zero. Thus $B_{l ijk}$ is homogeneous of degree $-1$.

Euler's formula applied to the second derivatives of $G$ gives $y^iB^m{}_{ijk}=0$. Complete symmetry of the lowered form supplies radial vanishing in the first slot as well. At a fixed $y\in\Sigma$, every $v\in V$ has the unique decomposition
\[
v=g_y(y,v)y+v^\top,\qquad v^\top\in T_y\Sigma.
\]
Radial vanishing thus makes the value of the full quartic equal to its value on the four tangential components. Its angular restriction determines it at $y$. Positive homogeneity determines it along the entire positive ray through $y$. Since $g_y$ is positive definite, raising or lowering an index preserves vanishing. This proves the last assertion.
\end{proof}

A smooth symmetric quartic on $(\Sigma,h)$ is a section of $\Sym^4T^*\Sigma$. Such a quartic $B$ is called \emph{twisted Codazzi} if, in any local tangent frame,
\begin{equation}\label{eq:twistedCodazzi}
(\nabla^h_pB)_{lijk}-(\nabla^h_lB)_{pijk}
=C_{li}{}^mB_{pmjk}-C_{pi}{}^mB_{lmjk}.
\end{equation}
Here $\nabla^h$ is the Levi--Civita connection and $C_{li}{}^m=h^{mq}C_{liq}$. All component indices range from $1$ to $n-1$ in this angular formula, and repeated indices are summed. Interchanging $i$ and $j$ gives
\begin{equation}\label{eq:angularCompatibility}
\begin{split}
0={}&C_{pi}{}^mB_{lmjk}-C_{li}{}^mB_{pmjk}\\
&-C_{pj}{}^mB_{lmik}+C_{lj}{}^mB_{pmik}.
\end{split}
\end{equation}
Indeed, the left-hand side of~\eqref{eq:twistedCodazzi} is unchanged by this interchange. Subtracting the right-hand sides gives~\eqref{eq:angularCompatibility}.

\begin{proposition}[Angular quartic rigidity]\label{prop:angular-rigidity}
Let $F$ be a smooth strongly convex Minkowski norm on an $n$-dimensional vector space $V$, regular on $V^\circ$, where $n\geq4$. Put $m=n-1$ and let $(\Sigma,h)$ be its indicatrix. Suppose that $B$ is a nonzero smooth symmetric quartic on $\Sigma$ satisfying~\eqref{eq:twistedCodazzi}.

At $y\in\Sigma$, define the insertion map, its radical, and its rank by
\[
\iota_y(v)=B_y(v,\cdot,\cdot,\cdot),\qquad
\rad B_y=\ker\iota_y,\qquad
\rank B_y=\rank\iota_y.
\]
Thus $\iota_y:T_y\Sigma\to\Sym^3T_y^*\Sigma$ and $\rank B_y=m-\dim\rad B_y$. Set
\[
q=\max_{y\in\Sigma}\rank B_y,\qquad
U=\{y\in\Sigma:\rank B_y=q\}.
\]
Then the following assertions hold.
\begin{enumerate}[label=\textup{(\roman*)}]
\item If $0<q<m$, then $q=m-1$ is even and $C$ vanishes on $U$.
\item If $q=m$, there is a point $p\in U$ such that
\[
\{A=A^*:[A,B_{uv}]=0\text{ for all }u,v\in T_p\Sigma\}=\R\Id,
\qquad h(B_{uv}w,z)=B(u,v,w,z).
\]
At this point there is a real number $s$ for which
\[
[C_X,C_Y]=sR_0(X,Y),\qquad X,Y\in T_p\Sigma.
\]
\end{enumerate}
Here $A^*$ denotes the $h_p$-adjoint. The displayed set in \textup{(ii)} is the self-adjoint common commutant of the contractions $B_{uv}$. As before, $C_XY=C^\sharp(X,Y)$ and $R_0(X,Y)Z=h(Y,Z)X-h(X,Z)Y$.

For a Landsberg metric, the angular restriction of $B_{l ijk}$ satisfies~\eqref{eq:twistedCodazzi}; hence the conclusions apply whenever this restriction is nonzero.
\end{proposition}

\begin{proof}
Part \textup{(i)} follows from Proposition~\ref{degenerate:conclusion}, and part \textup{(ii)} from Proposition~\ref{fullrank:generic} and Theorem~\ref{algebra:isotropy}, using~\eqref{eq:angularCompatibility}. The final assertion follows from Lemma~\ref{app:angular-equation}.
\end{proof}

We next apply Proposition~\ref{prop:angular-rigidity} to the Berwald quartic. The required curvature conclusion is intrinsic to $(\Sigma,h)$.

\begin{proposition}\label{prop:homogeneous-isotropy}
Suppose $K$ is transitive and the angular Berwald quartic is nonzero. Then the intrinsic metric $h$ has constant sectional curvature.
\end{proposition}
\begin{proof}
The maximum rank $q$ is positive, since the quartic is not identically zero. The maximum-rank set $U$ is nonempty and open: a nonzero maximal minor remains nonzero nearby, and no larger rank occurs by the definition of $q$.

If $q<m$, Proposition~\ref{prop:angular-rigidity} gives $C=0$ at some point $p\in U$. The angular Gauss formula gives $R^h_p=R_{0,p}$. If $q=m$, that proposition gives a point $p$ and scalar $s$ such that the same Gauss formula gives
\[
R^h_p=(1-s)R_{0,p}.
\]
In either case there is a real number $\kappa$ with $R^h_p=\kappa R_{0,p}$.

For an arbitrary $z\in\Sigma$, transitivity supplies $\phi\in K$ with $\phi(p)=z$. Isometries preserve the Levi--Civita connection and hence the curvature tensor. They also preserve $R_0$, which is constructed from $h$ alone. For $X,Y,Z\in T_p\Sigma$,
\[
R^h_z(d\phi X,d\phi Y)d\phi Z
=d\phi\bigl(R^h_p(X,Y)Z\bigr)
=\kappa R_{0,z}(d\phi X,d\phi Y)d\phi Z.
\]
Invertibility of $d\phi_p$ proves $R^h_z=\kappa R_{0,z}$ on the entire tangent space. Hence $h$ has constant sectional curvature $\kappa$.
\end{proof}

\begin{lemma}\label{lem:curvature-normalization}
Let $F$ be a smooth strongly convex Minkowski norm, regular on the entire punctured vector space, and suppose that its indicatrix has dimension $m\geq2$. If the indicatrix metric has constant sectional curvature $\kappa$, then $\kappa=1$.
\end{lemma}
\begin{proof}
Write $C_XY=C^\sharp(X,Y)$ and
\[
R_0(X,Y)Z=h(Y,Z)X-h(X,Z)Y.
\]
Our curvature convention and the Hessian cone calculation give
\begin{equation}\label{eq:constantGauss}
[C_X,C_Y]=(1-\kappa)R_0(X,Y).
\end{equation}
Let $\tau=\tr_h C$, a one-form. In fixed affine coordinates on $V$, the determinant derivative identity gives
\[
\partial_{y^k}\log\det g
=g^{ij}\partial_{y^k}g_{ij}=2g^{ij}C_{ijk}.
\]
On the indicatrix the radial term of this trace vanishes. Therefore
\begin{equation}\label{eq:tracePotentialHomogeneous}
\tau=2\dd f,\qquad f=\tfrac14\log\det(g_{ij})\big|_\Sigma.
\end{equation}
Changing affine coordinates adds only a constant to $f$. The potential is thus a globally defined smooth function after one choice of coordinates.

Contract~\eqref{eq:constantGauss} twice in an $h$-orthonormal basis $e_1,\ldots,e_m$. The contraction of the first composition is $|\tau|^2$, and that of the second is $|C|^2$, so
\begin{equation}\label{eq:traceConstantGauss}
m(m-1)(1-\kappa)=|\tau|_h^2-|C|_h^2.
\end{equation}
For example, the first contraction is
\[
\sum_{i,j}h(C_{e_i}C_{e_j}e_j,e_i)
=\left|\sum_i C_{e_i}e_i\right|^2=|\tau|^2,
\]
whereas the second is $\sum_{i,j}|C_{e_i}e_j|^2=|C|^2$.
Compactness of $\Sigma$ gives an extremum of $f$. At any such extremum, $\tau=0$, and~\eqref{eq:traceConstantGauss} implies $\kappa\geq1$.

Assume for a contradiction that $\kappa>1$, and choose a real number $t$ with
\[
t^2=\frac{\kappa}{\kappa-1}.
\]
Consider the connection on $T\Sigma$
\[
\nabla^t_XY=\nabla^h_XY+tC_XY.
\]
It is torsion-free because $C_XY=C_YX$. Its curvature is
\[
R^t(X,Y)=R^h(X,Y)
+t\{(\nabla^h_XC)_Y-(\nabla^h_YC)_X\}
+t^2[C_X,C_Y].
\]
The middle term vanishes by the ordinary Codazzi identity for the Cartan cubic. Using~\eqref{eq:constantGauss}, we obtain
\[
R^t=\{\kappa+t^2(1-\kappa)\}R_0=0.
\]
Thus $\nabla^t$ is a flat torsion-free connection on the entire indicatrix.

The radial parametrization identifies $\Sigma$ with $S^m$. Since $m\geq2$, it is simply connected. Flatness and simple connectivity make parallel transport independent of the path. A nonzero covector at one point therefore extends to a global $\nabla^t$-parallel one-form $\alpha$. Transport is an isomorphism, so $\alpha$ is nowhere zero. Torsion-freeness gives
\[
\dd\alpha(X,Y)=(\nabla^t_X\alpha)(Y)-(\nabla^t_Y\alpha)(X)=0.
\]
Every closed one-form on the simply connected sphere is exact, say $\alpha=\dd u$. The smooth function $u$ has an extremum on compact $\Sigma$, where $\dd u=0$. This contradicts the nonvanishing of $\alpha$. Hence $\kappa>1$ is impossible and $\kappa=1$.
\end{proof}

\begin{corollary}\label{cor:homogeneous-linearity}
On a connected component of the base with transitive indicatrices, either the Berwald curvature vanishes identically or every nonlinear curvature field $R(u,v)$ is linear in the original linear fibre coordinates.
\end{corollary}
\begin{proof}
If the Berwald curvature vanishes identically, the first alternative holds by definition. Otherwise select a fibre where it is nonzero. Lemma~\ref{lem:angularnonzero} supplies a nonzero angular quartic there. Proposition~\ref{prop:homogeneous-isotropy} and Lemma~\ref{lem:curvature-normalization} show that this indicatrix has constant sectional curvature one. Nonlinear parallel translation identifies it isometrically with all indicatrices of the component. Thus all of them have curvature one.

On any such fibre the fields $R(u,v)$ satisfy the smoothness, homogeneity, Killing, tangency, and fixed-label Bianchi assumptions of Theorem~\ref{round:linearity}. That theorem implies that each $R(u,v)$ is linear in the original fibre coordinates. This proves linearity of the nonlinear curvature.
\end{proof}

\section{Curvature fields on round indicatrices}\label{sec:round}

We now consider an indicatrix whose metric has constant sectional curvature one.

\begin{theorem}\label{round:linearity}\label{prop:round-linearity}
Let $V$ have dimension $n\geq3$, and let $E$ be a smooth, strongly convex,
positively two-homogeneous energy on $V^\circ=V\setminus\{0\}$. Set
$g=D^2E$, $C=\frac12Dg$, $\Sigma=\{2E=1\}$, and $h=g|_{T\Sigma}$.
Suppose that $h$ has constant sectional curvature one. Let
$R(u,v)$ be a family of smooth vector fields on $V^\circ$, alternating and
bilinear in the fixed labels $u,v\in V$, such that each field is positively
one-homogeneous, tangent to the levels of $E$, and Killing for $g$.
Suppose also that
\begin{equation}\label{round:bianchi}
 D_wR(u,v)+D_uR(v,w)+D_vR(w,u)=0
 \qquad(u,v,w\in V).
\end{equation}
Then every $R(u,v)$ is linear in the original variable $y\in V$.
\end{theorem}

The curvature fields constructed earlier satisfy
\eqref{round:bianchi} by \eqref{eq:VBianchi}. Throughout this section the
labels $u,v,w$ remain fixed vectors in the original vector space $V$. On the developed Euclidean space
introduced below, we write $\partial$ for ordinary differentiation, reserving
$D$ for ordinary differentiation on $V$.

\begin{definition}[Insertion rank and active axes]\label{def:cubic-rank}
For a symmetric cubic $T$ on a Euclidean space $W$, its insertion map is
$\iota_T:W\to\Sym^2W^*$, $\xi\mapsto T(\xi,\cdot,\cdot)$, and its
insertion rank is $\rank\iota_T$. A unit vector $e$ is an active axis if
$\lambda=T(e,e,e)\ne0$ and
\[
T(e,u,v)=\lambda\ip{e}{u}\ip{e}{v}\qquad(u,v\in W).
\]
The associated active line is $\R e$. Replacing $e$ by $-e$ reverses the
coefficient $\lambda$ and leaves the cubic summand $\lambda(e^*)^3$
unchanged. For a unit vector $e$, write $P_eu=\ip{u}{e}e$ for orthogonal
projection onto this line. Lemma~\ref{lem:round-cubic} shows that the commuting
insertion operators used below admit an orthonormal collection of active axes
representing the cubic.
\end{definition}

Write $r=\sqrt{2E(y)}$ and $\theta=y/r$. By homogeneity,
\[
 g_y(y,\cdot)=dE_y,\qquad g_y(y,y)=r^2,\qquad g_{r\theta}=g_\theta.
\]
Consequently the cone splitting \eqref{eq:cone} becomes
$g=dr^2+r^2h$. The radial parametrization identifies $\Sigma$ with the ordinary
sphere $S^{n-1}$. In particular it is compact and simply connected. The complete
simply connected space form theorem supplies an isometry
$\iota:(\Sigma,h)\longrightarrow S(W)$, where $W$ is an $n$-dimensional
Euclidean space and $S(W)$ is its unit sphere. Its radial extension
\begin{equation}\label{round:development}
 \Phi(r\theta)=r\iota(\theta)
 \colon V^\circ\longrightarrow W^\circ
\end{equation}
is a smooth positively one-homogeneous diffeomorphism and an isometry of the
cone metrics. Put $\Psi=\Phi^{-1}$, $Q_y=D\Phi_y$, and
$Q(z)=Q_{\Psi(z)}$. Thus $Q(z):V\longrightarrow W$ is invertible.
The map $\Phi$ is called a Euclidean development of the Hessian cone. Its homogeneity
implies that $Q$ has degree zero. The derivative $Q_y$ is taken in the
original linear coordinates on $V$; no linearity of $\Phi$ is assumed.

\begin{lemma}[Developed cubic]\label{lem:round-cubic}
There is a smooth $\End(W)$-valued one-form $A$ on $W^\circ$ such that
\begin{equation}\label{round:Q-equation}
 \partial_\xi Q=A(\xi)Q.
\end{equation}
Its insertions are self-adjoint, commute with one another, and satisfy
\begin{equation}\label{round:A-properties}
 A(\xi)\eta=A(\eta)\xi,\qquad A_z(z)=0,\qquad dA=0.
\end{equation}
The cubic
$\widehat C_z(\xi,\eta,\zeta)=\ip{A_z(\xi)\eta}{\zeta}$ is symmetric,
and $\partial\widehat C$ is symmetric in all four indices. At each point it
has an orthogonal decomposition
\begin{equation}\label{round:diagonal}
 \widehat C=\sum_{i=1}^q\lambda_i(e_i^*)^3,
 \qquad \lambda_i\ne0,\quad e_i\perp z,
\end{equation}
with $e_i$ orthonormal. Here $e_i^*=\ip{e_i}{\cdot}$ and $q$ is the insertion
rank, namely the dimension of the image of $\xi\mapsto A_z(\xi)$.
\end{lemma}

\begin{proof}
In fixed linear coordinates on $V$ the Levi--Civita connection of the Hessian metric
is $D+C^\sharp$, where
\[
 C^\sharp_y(u)v=g_y^{-1}C_y(u,v,\cdot),\qquad
 C^\sharp=\tfrac12g^{-1}dg.
\]
The isometry equation is $g=Q_y^TQ_y$. Naturality of the Levi--Civita
connection under $\Phi$ gives
$(D_uQ)_yv=Q_yC^\sharp_y(u)v$. Differentiating $g^{-1}dg$ yields
\[
 dC^\sharp(u,v)=-2[C^\sharp(u),C^\sharp(v)].
\]
The metric $g$ is flat, since it is isometric to the Euclidean punctured space.
Its curvature equation is
$dC^\sharp(u,v)+[C^\sharp(u),C^\sharp(v)]=0$. Subtracting these equations shows
that $dC^\sharp=0$ and all insertion commutators vanish.

For $y=\Psi(z)$ define
\[
 A_z(\xi)=Q_yC^\sharp_y(Q_y^{-1}\xi)Q_y^{-1}.
\]
The chain rule gives \eqref{round:Q-equation}. Self-adjointness follows by
transporting $g$ to the Euclidean inner product; cubic symmetry gives
$A(\xi)\eta=A(\eta)\xi$. Euler's identity for $\Phi$ is $Q_yy=\Phi(y)$.
Together with $C_y(y,\cdot,\cdot)=0$, it implies $A_z(z)=0$.
Finally $A=(dQ)Q^{-1}$ is a right Maurer--Cartan form. Hence
$dA(\xi,\eta)=[A(\xi),A(\eta)]=0$. This closedness, combined with the
symmetry of $\widehat C$, gives total symmetry of $\partial\widehat C$.

Simultaneously diagonalize the commuting self-adjoint operators $A(\xi)$.
For two different common unit eigenvectors $e_i,e_j$, symmetry says that
$A(e_i)e_j=A(e_j)e_i$. The first vector lies in $\R e_j$, the second in
$\R e_i$, so both vanish. Only pure cubic coefficients remain, giving
\eqref{round:diagonal}. Radial nullity implies
$\lambda_i\ip{e_i}{z}=0$ for every retained term. The nonzero diagonal
insertions are independent, proving the description of $q$.
\end{proof}

Choose linear coordinates on $V$ and an orthonormal frame on $W$. Define
\[
 \varphi(z)=\log|\det Q(z)|.
\]
Since $g=Q^TQ$, the function $f=\tfrac14\log\det(g_{ij})$ satisfies
$\varphi=2f\circ\Psi$. Changing the chosen volume forms adds a constant
to $\varphi$. The trace one-form of the developed Cartan cubic is
$\widehat\tau(\xi)=\tr A(\xi)$; thus $d\varphi=\widehat\tau$.
Since $Q$ is homogeneous of degree zero, so is $\varphi$. By
\eqref{round:Q-equation} and \eqref{round:diagonal},
\begin{equation}\label{round:gradient}
 d\varphi(\xi)=\tr A(\xi),\qquad
 G:=\operatorname{grad}\varphi=\sum_i\lambda_i e_i,\qquad
 \norm{G}^2=\norm{\widehat C}^2.
\end{equation}
Here the cubic norm is the Euclidean tensor norm. In particular $G$ is tangent
to every sphere, and its zeros are exactly the zeros of $\widehat C$.

\begin{definition}[Connection associated with a matrix-valued one-form]\label{def:flat-pencil}
For a real parameter $s$, the notation $d+sA$ denotes the connection on
the trivial bundle $W^\circ\times W$ given by
$(d+sA)_\xi v=\partial_\xi v+sA(\xi)v$. It is flat if its curvature
$s\,dA(\xi,\eta)+s^2[A(\xi),A(\eta)]$ is zero for all $\xi,\eta$.
A matrix whose columns form a parallel frame therefore satisfies
$\partial_\xi H=-sA(\xi)H$.
\end{definition}

For every fixed real number $s$, the connection $d+sA$ is flat:
its curvature is $s\,dA(\xi,\eta)+s^2[A(\xi),A(\eta)]=0$.
Since $W^\circ$ is simply connected, there is a global smooth invertible
matrix $H_s$ satisfying
\begin{equation}\label{round:pencil}
 H_s(z_*)=\Id,\qquad \partial_\xi H_s=-sA(\xi)H_s.
\end{equation}
For each fixed $s$, compactness of $S(W)$ bounds $H_s$ there. The equation
$A_z(z)=0$ also makes $H_s$ constant on positive rays.

The pushforward $\Phi_*R(u,v)$ is a positively one-homogeneous Euclidean
Killing field on $W^\circ$. Such a field is linear and skew-adjoint.
Indeed, differentiating the equations
$\partial_iX_j+\partial_jX_i=0$ and permuting the derivative indices gives
$\partial_i\partial_jX_k=0$. Connectedness makes $X$ affine, and positive
homogeneity removes its constant term. We therefore obtain a fixed alternating
bilinear map $\beta:V\times V\longrightarrow\so(W)$ with
\begin{equation}\label{round:beta}
 R(u,v)(y)=Q_y^{-1}\beta(u,v)\Phi(y).
\end{equation}
Differentiating this formula along a constant vector $w\in V$, multiplying by
$Q_y$, and using \eqref{round:Q-equation} transforms the vertical Bianchi
identity into
\begin{equation}\label{round:developed-bianchi}
 \sum_{\cyc(u,v,w)}\beta(u,v)Q_yw
 =\sum_{\cyc(u,v,w)}A_z(Q_yw)\beta(u,v)z,
 \qquad z=\Phi(y).
\end{equation}

The insertion rank need not be locally constant. The following lemma gives
smooth continuations of the nonzero eigendirections and estimates the
remaining cubic components.

\begin{lemma}[Selected axes and the residual term]\label{lem:round-residual}
Suppose that $\widehat C$ has rank $q>0$ at $z_0$, with ordered active axes
$e_1,\ldots,e_q$ and coefficients $\lambda_i$. On a neighborhood of $z_0$
these axes have smooth unit continuations and
\begin{equation}\label{round:residual-decomposition}
 \widehat C=\sum_{i=1}^q\lambda_i(e_i^*)^3+\mathcal E,
 \qquad
 \mathcal E=\mathcal E(P^\perp\cdot,P^\perp\cdot,P^\perp\cdot),
\end{equation}
where $P^\perp$ projects off the selected axes. The selected coefficients
remain nonzero on a smaller neighborhood. In the multilinear operator
norm,
\begin{equation}\label{round:residual-estimate}
 \norm{\partial_{e_i}\mathcal E}
 \leq
 \bigl(\norm{\partial e_i}
       +3\norm{\partial_{e_i}P^\perp}\bigr)\norm{\mathcal E}.
\end{equation}
The same statement holds when one selects any subset of the active axes.
\end{lemma}

\begin{proof}
At $z_0$ fix the vector $w_0=\sum_{i=1}^q i e_i/\lambda_i$.
The symmetric matrix $A_{z_0}(w_0)$ has simple nonzero eigenvalues
$1,\ldots,q$, with zero on the remaining subspace. Its simple eigenlines
continue smoothly nearby, with signs fixed by positive overlap with their
initial values (the standard simple-eigenvalue perturbation argument; see~\cite{Kato}). Every insertion commutes with $A(w_0)$, so it preserves each
selected eigenline. Cubic symmetry then annihilates every mixed coefficient
involving a selected line. This proves \eqref{round:residual-decomposition}
without any restriction on the residual rank. For a subset of axes, use
distinct nonzero weights on those axes and zero on the others.

At a reference point choose constant vectors $X,Y,Z$ perpendicular to the
selected axes there. Differentiate
\[
 \widehat C(e_i,Y,Z)
 =\lambda_i\ip{e_i}{Y}\ip{e_i}{Z}
\]
in the direction $X$. The right-hand derivative vanishes at the reference point.
Since $\partial\widehat C$ is totally symmetric,
\[
 (\partial_{e_i}\widehat C)(X,Y,Z)
 =-\widehat C(\partial_Xe_i,Y,Z)
 =-\mathcal E(\partial_Xe_i,Y,Z).
\]
The derivative of each selected pure cubic on three perpendicular inputs is
zero. Thus the restriction of $\partial_{e_i}\mathcal E$ to such inputs is
bounded by $\norm{\partial e_i}\norm{\mathcal E}$. To control arbitrary
inputs, differentiate the support identity in
\eqref{round:residual-decomposition}. The three differentiated projections
contribute at most
$3\norm{\partial_{e_i}P^\perp}\norm{\mathcal E}$.
This gives \eqref{round:residual-estimate}. Here $\norm{\partial e_i}$ means
the norm of the full derivative $X\mapsto\partial_Xe_i$.
\end{proof}

Active axes are locally unique up to sign. Indeed, suppose
$e,f$ are two unit axes with nonzero coefficients $\lambda,\mu$. Symmetry
gives
$\lambda\ip{e}{f}e=A(e)f=A(f)e=\mu\ip{e}{f}f$.
Thus two distinct active lines are orthogonal. In particular two smoothly
continued axes close to the same initial axis agree after their signs are
fixed. This uniqueness ensures compatibility of the local continuations of the ordered eigenvectors used below.

\begin{definition}[Regular point and regular value]\label{def:regular-level}
A point $z\in S(W)$ is regular for $\varphi|_{S(W)}$ if its differential on
$T_zS(W)$ is nonzero, equivalently if $G(z)\ne0$. A real number $c$ is a
regular value if every point of the level $\{\varphi=c\}\cap S(W)$ is regular.
A nonempty regular level is a smooth hypersurface of $S(W)$ and is compact
because it is closed in the unit sphere.
\end{definition}

If $\widehat C$ vanishes identically, then \eqref{round:Q-equation} makes
$Q$ constant, and \eqref{round:beta} proves the conclusion.
Assume henceforth that $\widehat C\not\equiv0$. Equation
\eqref{round:gradient} shows that $\varphi|_{S(W)}$ is nonconstant. Sard's theorem
therefore supplies a nonempty compact regular level
$L=\{z\in S(W):\varphi(z)=c\}$. Let $q\geq1$ be the least insertion rank on
$L$, and form the framed set
\[
 \mathcal K=\left\{(z,e_1,\ldots,e_q):
 \begin{array}{l}
 z\in L,\ e_i\text{ orthonormal},\\
 \widehat C_z=\sum_{i=1}^q\widehat C_z(e_i,e_i,e_i)(e_i^*)^3
 \end{array}\right\}.
\]
The set $\mathcal K$ consists of points of $L$ together with ordered
orthonormal active axes representing the entire cubic. It is nonempty.
Its defining identities are closed in
$L\times S(W)^q$, so $\mathcal K$ is compact. No selected coefficient can
vanish there, since that would produce rank less than $q$ on $L$.
Consequently all $|\lambda_i|$ have a common positive lower bound on
$\mathcal K$.

\begin{lemma}\label{lem:round-rank-one}
The minimum rank $q$ on $L$ equals one.
\end{lemma}

\begin{proof}
Suppose $q\geq2$ and fix two distinct ordered indices $i,j$. In a neighborhood
of a point of $\mathcal K$, use Lemma~\ref{lem:round-residual} and set
\[
 V=\frac{e_i}{\lambda_i}-\frac{e_j}{\lambda_j}.
\]
This vector is tangent to the unit sphere. Also $d\varphi(V)=1-1=0$, since
the residual trace is perpendicular to the selected axes. These assertions
hold on the neighborhood even where $\mathcal E$ does not vanish.
The residual estimate gives, along a local integral curve of $V$,
\[
 \norm{\tfrac d{dt}\mathcal E}\leq c\norm{\mathcal E}
\]
with a locally bounded coefficient. Gronwall's inequality in either time
direction shows that $\mathcal E=0$ is preserved. Therefore a solution
starting in $\mathcal K$, with its continued ordered axes, stays in
$\mathcal K$ for its local lifetime.

We next prove that these curves are complete. The argument does not require
a smooth-manifold structure on $\mathcal K$.
Near each point of $\mathcal K$, all sufficiently close framed points of
$\mathcal K$ lie on the graph of the smooth eigenline continuations used in
the preceding lemma. The lifted
local vector is
\[
 \bigl(V,\partial_Ve_1,\ldots,\partial_Ve_q\bigr).
\]
On overlapping neighborhoods the continued eigenlines agree near a common
framed point by uniqueness of an active axis with its prescribed sign.
Thus the local solutions agree and are unique. A finite cover of the compact
set gives a uniform speed bound for the lifted curves. At a finite maximal
endpoint, bounded speed makes the curve Cauchy in the ambient product; its
limit belongs to $\mathcal K$ and a local solution there continues it.
The lifted solution is therefore defined for all $t\geq0$. In particular,
there is a finite constant $b$ with $\norm{e_j'(t)}\leq b$ along it.

Choose one fixed $s>2b$ in \eqref{round:pencil}, and choose $v_*$ so that
$h(t)=H_s(z(t))v_*$ satisfies $h(0)=e_j(0)$. This is possible because
$H_s(z(0))$ is invertible. Let $P_i,P_j$ be the orthogonal projections
onto the selected axes. Since $\mathcal E=0$ on the lifted curve,
\[
 A(V)=P_i-P_j,\qquad h'=s(P_j-P_i)h.
\]
Put $f=\ip{h}{e_j}$, $p=\ip{h}{e_i}$, and $U=2f^2-\norm{h}^2$.
Direct differentiation gives
\[
 U'=2s(f^2+p^2)+4f\ip{h}{e_j'}.
\]
Whenever $U\geq0$, one has
$\norm{h-fe_j}\leq|f|$ and $e_j'\perp e_j$. Therefore
\[
 U'\geq(2s-4b)f^2\geq(2s-4b)U.
\]
Since $U(0)=1$, a continuity argument shows that $U$ remains positive, and integration
gives $U(t)\geq\exp((2s-4b)t)$ for all $t\geq0$.
On the other hand, $H_s$ is bounded on the unit sphere, so $h$ and hence
$U$ are bounded. This contradiction proves $q=1$.
\end{proof}

Start at a rank-one point $z(0)\in L$ and follow
\begin{equation}\label{round:gradient-arc}
 z'=\frac{G}{\norm{G}^2}
\end{equation}
on its maximal forward interval $[0,T_*)$ in the regular set of $\varphi$ on
the sphere. The rank remains one along this interval. Indeed, near a regular
rank-one point select its active axis, writing
$\widehat C=\lambda e^3+\mathcal E$. Then
$G=\lambda e+\tr\mathcal E$. The trace, regarded as a vector, is perpendicular
to $e$. Its norm is bounded by $\norm{\mathcal E}$ times a constant depending
only on the dimension. Lemma~\ref{lem:round-residual} controls
$\partial_e\mathcal E$; smoothness controls the derivative in the direction of
$\tr\mathcal E$. Since $\norm{G}$ is locally bounded below,
\eqref{round:gradient-arc} gives
$\norm{\mathcal E'}\leq c\norm{\mathcal E}$.
Thus $\mathcal E=0$ is locally preserved. Rank at most one is a closed
condition, and $G\ne0$ excludes rank zero. An open-and-closed argument on
the maximal interval proves the assertion.

We have $\varphi(z(t))=\varphi(z(0))+t$, so $T_*<\infty$ because $\varphi$ is bounded
on the compact sphere. Every cluster point at the forward endpoint is
critical. Otherwise the smooth vector field \eqref{round:gradient-arc}
has a uniform short existence interval near a regular cluster point, extending
the same trajectory past $T_*$ by uniqueness. At a critical cluster point,
$\widehat C=0$ by \eqref{round:gradient}.

Along the arc choose the smooth orientation
$e=G/\norm{G}$, $\lambda=\norm{G}>0$. Then
\begin{equation}\label{round:arc-Q}
 z'=e/\lambda,\qquad A(z')=P_e,\qquad Q'=P_eQ.
\end{equation}
For every triple of fixed labels $u,v,w\in V$, set
\[
 \mathcal D(u,v,w)(t)=\sum_{\cyc(u,v,w)}\beta(u,v)Q(t)w.
\]
Equation \eqref{round:developed-bianchi} shows that $\mathcal D(u,v,w)$ is a multiple of $e$.
Its derivative is $\sum_{\cyc}\beta(u,v)P_eQ(t)w$, which is perpendicular
to $e$ because each $\beta(u,v)$ is skew-adjoint. Hence
$\norm{\mathcal D(u,v,w)}^2$ is constant. Along a subsequence approaching a critical
cluster point, the right-hand side of
\eqref{round:developed-bianchi} tends to zero; $Q$ and the other factors
remain bounded on the sphere. The constant is therefore zero:
\begin{equation}\label{round:algebraic-bianchi}
 \sum_{\cyc(u,v,w)}\beta(u,v)Q(t)w=0.
\end{equation}

Identify $\Lambda^2W$ with $\so(W)$ by
\[
 u\wedge v\longmapsto
 \bigl(w\longmapsto\ip{v}{w}u-\ip{u}{w}v\bigr).
\]
Here $\so(W)$ is the space of skew-adjoint endomorphisms of $W$. The
Euclidean inner product on $W$ induces the inner product on $\Lambda^2W$
defined on decomposable bivectors by
\[
\ip{u\wedge v}{w\wedge z}
=\ip{u}{w}\ip{v}{z}-\ip{u}{z}\ip{v}{w}.
\]
Adjoints and orthogonal projections on $\Lambda^2W$ are taken with respect
to this inner product.
Using this identification, define the relative curvature operator
\begin{equation}\label{round:relative}
 \mathcal R_t(u\wedge v)
 =\beta(Q(t)^{-1}u,Q(t)^{-1}v).
\end{equation}
This operator acts on $\Lambda^2W$. Equation
\eqref{round:algebraic-bianchi} is its first algebraic Bianchi identity.
Skewness in each pair and that identity imply interchange of the pairs,
so $\mathcal R_t$ is self-adjoint on $\Lambda^2W$. For completeness, if
$T_{ijkl}$ denotes the associated four-tensor and
$x=T_{ijkl}$, $y=T_{ikjl}$, $z=T_{iljk}$,
$p=T_{jkil}$, $q=T_{jlik}$, $r=T_{klij}$, the four Bianchi identities give
\[
 p=y-x,\qquad q=x+z,\qquad r=y-z,\qquad q=p+r.
\]
They imply $x=r$, which is the required pair interchange.

Let $\Pi_t$ be the orthogonal projection onto
$e(t)\wedge W\subset\Lambda^2W$. Explicitly,
\[
 \Pi_t(u\wedge v)=(P_eu)\wedge v+u\wedge(P_ev).
\]
The decomposition
$\Lambda^2W=(e\wedge e^\perp)\oplus\Lambda^2e^\perp$ verifies that this
formula defines an orthogonal projection. Differentiating
\eqref{round:relative} and using \eqref{round:arc-Q} gives
\begin{equation}\label{round:operator-ODE}
 \mathcal R'=-\mathcal R\Pi,\qquad [\mathcal R,\Pi]=0.
\end{equation}
The commutation follows because both $\mathcal R$ and $\mathcal R'$ are
self-adjoint.

\begin{lemma}[Projection rigidity on the curvature image]
\label{lem:round-projection}
The subspace $I=\im\mathcal R_t$ is independent of $t$, and
$\Pi_t|_I$ is constant.
\end{lemma}

\begin{proof}
Solve $\mathcal U'=-\mathcal U\Pi$, $\mathcal U(0)=\Id$.
Its fundamental matrix is invertible, and uniqueness gives
$\mathcal R_t=\mathcal R_0\mathcal U_t$. Hence the image is fixed.
Self-adjointness gives $\ker\mathcal R_t=I^\perp$, and commutation in
\eqref{round:operator-ODE} makes $I$ and $I^\perp$ invariant under $\Pi_t$.
They are therefore invariant under its derivatives as well.

Differentiating $[\mathcal R,\Pi]=0$ and using
$[\mathcal R',\Pi]=0$ yields $[\mathcal R,\Pi']=0$.
Differentiating again gives
\[
 [\mathcal R,\Pi'']=-[\mathcal R',\Pi']
 =\mathcal R[\Pi,\Pi'].
\]
At a fixed time decompose into eigenspaces of the self-adjoint
$\mathcal R$. The left side has zero diagonal blocks. The right side
commutes with $\mathcal R$, because both $\Pi$ and $\Pi'$ do, and is
therefore block diagonal. Both sides must vanish. On $I$ the operator
$\mathcal R$ is invertible, so $[\Pi,\Pi']=0$ there.
Differentiating $\Pi^2=\Pi$ gives
$\Pi\Pi'+\Pi'\Pi=\Pi'$. Commutation on $I$ reduces this to
$2\Pi\Pi'=\Pi'$. Multiplying by $\Pi$ first shows $\Pi\Pi'=0$ on $I$,
and then $\Pi'=0$ on $I$.
\end{proof}

The active lines cannot all be parallel. Otherwise
\eqref{round:arc-Q} would place a nonconstant spherical arc in an affine
line, whose intersection with the unit sphere has at most two points.
Choose two times with nonparallel axes, let $W_1$ be their span, and put
$W_2=W_1^\perp$. Both summands are nonzero, since $n\geq3$.
By Lemma~\ref{lem:round-projection}, the fixed image of $\Pi|_I$ lies in
the intersection of the two spaces $e\wedge W$, namely $\Lambda^2W_1$.
Its fixed kernel lies in the intersection of their orthogonal complements,
namely $\Lambda^2W_2$. Thus
\[
 I\subset\Lambda^2W_1\oplus\Lambda^2W_2.
\]
The mixed subspace $W_1\wedge W_2$ is orthogonal to $I$, so
\begin{equation}\label{round:mixed-zero}
 \mathcal R_0(W_1\wedge W_2)=0.
\end{equation}

\begin{proof}[Proof of Theorem~\ref{round:linearity}]
Put $V_i=Q(0)^{-1}W_i$, so $V=V_1\oplus V_2$. By
\eqref{round:relative} and \eqref{round:mixed-zero},
$\beta(V_1,V_2)=0$. Formula \eqref{round:beta} therefore shows that every
field with labels in different summands vanishes on the entire slit space $V^\circ$.
For $v_1,v_2\in V_1$ and $w\in V_2$, the fixed-label identity
\eqref{round:bianchi} now gives $D_wR(v_1,v_2)=0$.

For fixed $v_1,v_2\in V_1$, the field $R(v_1,v_2)$ is smooth and positively
one-homogeneous on $V^\circ$, and $D_wR(v_1,v_2)=0$ for every $w\in V_2$.
Lemma~\ref{lem:quotient-origin} therefore shows that it is the restriction of
a linear map $V\to V$ that vanishes on $V_2$. Interchanging $V_1,V_2$ proves
the analogous statement for labels in $V_2$.
Bilinearity and alternation
therefore give linearity for arbitrary labels.
\end{proof}

\begin{remark}
Lemma~\ref{lem:quotient-origin} uses smoothness at every nonzero vector of
$V$; positive homogeneity alone does not imply linearity at the quotient origin.
The conclusion concerns the curvature family and follows from its fixed-label
Bianchi identity.
\end{remark}

\section{Proof of the main theorem}\label{sec:compact}

We first identify the $R$-quadratic condition supplied by the preceding fibrewise arguments, and then apply Crampin's compact rigidity theorem~\cite{Crampin2011,Crampin2024}.

Suppose every curvature field $R(u,v)$ is linear in the fibre variable. In a tangent chart,
\begin{equation}\label{eq:raw-linear}
R^m{}_{ij}(x,y)=A^m{}_{ijp}(x)y^p.
\end{equation}
The coefficients depend smoothly on $x$: for any smooth nonzero local section $s(x)$, they are given by $A^m{}_{ijp}(x)=\partial_{y^p}R^m{}_{ij}(x,s(x))$. The derivative is independent of the section because the curvature is linear on each fibre. Thus $A$ is a smooth section of $\Lambda^2T^*M\otimes\End(TM)$.

For a fixed vector $v\in T_xM$ write $\rho_v(y)=R(y,v)(y)$. Let $\partial_u$ denote ordinary fibre differentiation in a constant direction $u$. The product rule and the fixed-argument vertical Bianchi identity give
\begin{align*}
\partial_u\rho_v-\partial_v\rho_u
&=2R(u,v)+\partial_uR(y,v)-\partial_vR(y,u),\\
0&=(\partial_yR(u,v))(y)-\partial_uR(y,v)+\partial_vR(y,u).
\end{align*}
In the second line, the vector equal to the reference point $y$ is held fixed both as an argument and as a differentiation direction. Euler's identity for the one-homogeneous field $R(u,v)$ therefore yields
\begin{equation}\label{eq:rho-factor-three}
\partial_u\rho_v-\partial_v\rho_u=3R(u,v).
\end{equation}
Consequently the curvature fields are linear if and only if the Jacobi endomorphism is a quadratic polynomial in $y$. This is the $R$-quadratic condition of Definition~\ref{def:R-quadratic}.

The Cartan tensor retains the normalization used throughout this paper:
\[
C_{ijk}=\frac12\frac{\partial g_{ij}}{\partial y^k}.
\]
\begin{proposition}[Crampin~\cite{Crampin2011,Crampin2024}]\label{prop:compact-terminal}
Every regular smooth strongly convex $R$-quadratic Landsberg metric on a closed manifold of dimension at least three is Berwald.
\end{proposition}

This is the compact-base conclusion of \cite[Theorem~1]{Crampin2011}; see also \cite[Theorem~5]{Crampin2024}. By the preceding argument, linearity of all curvature fields $R(u,v)$ verifies its $R$-quadratic hypothesis.

\begin{proof}[Proof of Theorem~\ref{thm:main}]
Assume first that $M$ is connected. Nonlinear parallel translation is an isometry between indicatrix metrics. Hence transitivity of $\Isom(\Sigma_x,h_x)^0$ is independent of $x$.

Suppose first that the action is nontransitive. Theorem~\ref{thm:proper-linearity} shows that every curvature field is linear in the fibre variable. Thus $F$ is $R$-quadratic, and Proposition~\ref{prop:compact-terminal} implies that $F$ is Berwald.

Suppose next that the action is transitive. If $B$ vanishes identically, there is nothing to prove. Otherwise choose a fibre on which $B$ is nonzero. Proposition~\ref{prop:angular-rigidity} and the curvature normalization in Section~\ref{sec:homogeneous} show that its indicatrix is isometric to the unit sphere. By nonlinear parallel translation, the same is true of every indicatrix. Theorem~\ref{round:linearity} shows that every curvature field is linear in the original fibre coordinates. Again $F$ is $R$-quadratic, so Proposition~\ref{prop:compact-terminal} gives the Berwald property.

If $M$ is disconnected, apply the argument to each connected component. Each component is closed in $M$ and therefore compact. This proves the theorem.
\end{proof}

\begin{remark}
Compactness of the base enters only through Crampin's theorem. The preceding fibrewise arguments use compactness of each indicatrix, which holds for every regular strongly convex Minkowski norm independently of the base.
\end{remark}

\appendix
\section{Notation and homogeneity}\label{app:dictionary}

The differentiation conventions are summarized below. On a fixed tangent space, $D$ is the flat affine derivative. Horizontal Berwald derivatives are denoted by a vertical bar, and the Cartan tensor is consistently normalized by $C=\tfrac12Dg$ throughout the paper.

\begin{center}
\renewcommand{\arraystretch}{1.25}
\begin{tabular}{>{\raggedright\arraybackslash}p{0.16\textwidth}>{\raggedright\arraybackslash}p{0.33\textwidth}>{\raggedright\arraybackslash}p{0.40\textwidth}}
\toprule
Notation & Domain & Meaning\\
\midrule
$D$ & A fixed vector space $V=T_xM$ & Flat affine derivative in the fibrewise arguments and Appendix~\ref{app:angular-proof}.\\
$\nabla^g$ & $TV^\circ$ & Levi--Civita connection of the fibre Hessian metric; $D+C^\sharp$.\\
$\nabla^h$ & $T\Sigma$ & Levi--Civita connection of the indicatrix metric.\\
$\nablac$ & $T\Sigma$ & Induced centroaffine connection $\nabla^h-C^\sharp$.\\
$\nabla^B$ & Pullback of $TM$ to $T^\circ M$ & Berwald connection from the canonical spray.\\
$\partial$ & Developed Euclidean space $W$ & Flat derivative after the Euclidean development $\Phi$.\\
$\nabla^t$ & $T\Sigma$ & Auxiliary connection $\nabla^h+tC^\sharp$ used only in normalization.\\
$\dd+sA$ & Trivial bundle $W^\circ\times W$ & Family of flat connections used in Section~\ref{sec:round}.\\
\bottomrule
\end{tabular}
\end{center}

For each $s\in\R$, the connection $\dd+sA$ on the developed fibre has zero curvature. The symbol $R^g$ denotes the Riemannian curvature of the Hessian metric, $R^h$ the curvature of the indicatrix metric, and $R(u,v)$ the curvature fields of the nonlinear connection with fixed vector arguments. The horizontal Berwald curvature is $H^m{}_{pij}=-\partial_{y^p}R^m{}_{ij}$. The operator $\mathcal R_t$ of Section~\ref{sec:round} acts on $\Lambda^2W$.

The components of the principal tensors have the following degrees of homogeneity in fixed linear fibre coordinates:
\begin{center}
\renewcommand{\arraystretch}{1.15}
\begin{tabular}{lll}
\toprule
Object & Degree & Radial identity\\
\midrule
$E$ & $2$ & $\dd E_y(y)=2E(y)$\\
$g=D^2E$ & $0$ & $g_y(y,\cdot)=\dd E_y$\\
$C=\tfrac12Dg$ & $-1$ & $C_y(y,\cdot,\cdot)=0$\\
$G^m$ & $2$ & $y^i\partial_{y^i}G^m=2G^m$\\
$N^m_i$ & $1$ & $N^m_i y^i=2G^m$\\
$\Gamma^m_{ij}$ & $0$ & $\Gamma^m_{ij}y^j=N^m_i$\\
$B^m{}_{ijk}$ & $-1$ & $B^m{}_{ijk}y^i=0$\\
$R^m{}_{ij}$ & $1$ & $D(R(u,v))_y[y]=R(u,v)(y)$ for fixed arguments\\
$\rho_{x,y}(v)=R(y,v)(y)$ & $2$ & The first argument also varies\\
\bottomrule
\end{tabular}
\end{center}

All argument vectors are held fixed in these homogeneity statements. For example, $B_{ty}(v_1,v_2,v_3,v_4)=t^{-1}B_y(v_1,v_2,v_3,v_4)$ when the test vectors are held fixed. If those vectors are dilated as well, four additional factors of $t$ appear. The proofs use component homogeneity in the former sense.

For variable vector arguments $u(y)$ and $v(y)$, the product rule is
\[
D_w\{R(u(y),v(y))(y)\}
=(D_wR(u,v))(y)+R(D_wu,v)(y)+R(u,D_wv)(y).
\]
The vertical Bianchi identity applies to the first term. The additional argument derivatives account for the factor three in~\eqref{eq:rho-factor-three}.

The tensor rank used in the transitive case is the rank of
\[
v\longmapsto B(v,\cdot,\cdot,\cdot).
\]
In Section~\ref{sec:round}, insertion rank means the dimension of the image of $\xi\mapsto A(\xi)$ for a commuting cubic. Neither is the dimension of the pointwise span of curvature-field values. The space $\mathcal W\subset\Lambda^2V^*$ is the coefficient space of the family $R(u,v)$, regarded as vector fields on the entire punctured fibre.

\section{Regularity and dimension assumptions}\label{app:variants}

For a direct sum $V=U\oplus Z$, the quotient vector space $V/Z$ consists of equivalence classes modulo $Z$ and is identified with $U$. A function constant along the affine subspaces $u+Z$ descends to this quotient wherever it is defined. We use the following extension lemma in Section~\ref{sec:round}.

\begin{lemma}\label{lem:quotient-origin}
Let $V=U\oplus Z$ with both summands nonzero. Let $f:V\setminus0\to W$ be smooth, positively one-homogeneous, and satisfy $D_zf=0$ for every constant $z\in Z$. Then $f$ is the restriction of a linear map $V\to W$ that vanishes on $Z$.
\end{lemma}
\begin{proof}
For $u\ne0$, the affine slice $u+Z$ avoids zero and is connected. Vanishing derivatives make $f$ constant there. On $Z\setminus0$, its derivative in the radial direction vanishes, while homogeneity gives $D_zf(z)=f(z)$. Hence $f(z)=0$ on that set. Fix $z_*\ne0$ in $Z$ and put $\widetilde f(u)=f(u+z_*)$. It is smooth at $u=0$, has value zero there, and agrees with $f(u)$ for $u\ne0$. It is positively one-homogeneous, including at zero by continuity. Differentiability at zero gives
\[
\widetilde f(u)=\lim_{t\downarrow0}\frac{\widetilde f(tu)-\widetilde f(0)}t
=D\widetilde f(0)u.
\]
Thus $f(u+z)=D\widetilde f(0)u$, including $u=0$. This is the required linear map.
\end{proof}

Evaluation on the affine subspace $U+z_*$ converts smoothness at nonzero vectors into differentiability at the quotient origin. Without this extension property, positive homogeneity does not imply linearity. For example, the map $u\mapsto |u|$ is smooth on $U\setminus0$, is positively one-homogeneous, and is nonlinear. Its pullback to $\{u+z:u\ne0\}$ has all $Z$ derivatives zero, but it does not extend smoothly across points with $u=0$, $z\ne0$. This scalar example illustrates the extension property used in the lemma.

In the proper-orbit argument, the exterior-product symbol and the coefficient two-form space are defined in every dimension. The transitive case uses Proposition~\ref{prop:angular-rigidity}, which assumes $n\geq4$; its proof is given in Appendix~\ref{app:angular-proof}. The round-linearity theorem requires $n\geq3$, which ensures simple connectivity of the punctured space and of the indicatrix. The curvature-normalization argument requires angular dimension at least two for the same topological reason.

In dimension three, the nontransitive arguments of Sections~\ref{sec:shape} and~\ref{sec:proper} remain valid. More explicitly, Lemma~\ref{lem:good-symbol} makes
\[
\mathcal W\longrightarrow\Lambda^3V^*,\qquad
\omega\longmapsto\xi\wedge\omega
\]
injective for a suitable nonzero covector $\xi$, where $\mathcal W$ is the coefficient space defined in~\eqref{eq:Wdef}. Since $\dim\Lambda^3V^*=1$, one has $\dim\mathcal W\leq1$, and the curvature fields are consequently linear. In the transitive case, curvature isotropy at a point is automatic because the indicatrix is a surface. Transitivity makes the curvature scalar constant, Lemma~\ref{lem:curvature-normalization} normalizes it to one, and Theorem~\ref{prop:round-linearity} again gives linearity. Thus the closed three-dimensional Berwald conclusion is recovered without Proposition~\ref{prop:angular-rigidity}.

Mao--Zhou~\cite{MaoZhou} additionally prove the three-dimensional Laugwitz conjecture by a separate argument. The present curvature-field method does not establish that stronger conclusion: for a translation-invariant Minkowski metric, the curvature fields are already zero, so their linearity gives no further restriction on the Cartan cubic. In dimension two, the global development and flat-pencil construction used here are not justified by simple connectivity.

\section{Proof of angular quartic rigidity}\label{app:angular-proof}

This appendix proves Proposition~\ref{prop:angular-rigidity}. Fix a smooth strongly convex Minkowski norm on a vector space $V$ of dimension $n=m+1\geq4$, regular on all of $V^\circ$. We retain $E=F^2/2$, $\Sigma=\{F=1\}$, $h=g|_{T\Sigma}$, and the angular Cartan cubic $C$. We write $C_XY=C^\sharp(X,Y)$, so that $h(C_XY,Z)=C(X,Y,Z)$. Let $B$ be a nonzero smooth symmetric quartic on the entire indicatrix satisfying~\eqref{eq:twistedCodazzi}. Its insertion rank, radical, and maximum-rank set are those defined in Proposition~\ref{prop:angular-rigidity}.

For substitutions involving subbundles, we use the following vector forms of~\eqref{eq:twistedCodazzi} and~\eqref{eq:angularCompatibility}:
\begin{equation}\label{app:vector-TC}
\begin{split}
(\nablah_XB)(Y,Z,U,W)-(\nablah_YB)(X,Z,U,W)
={}&B(X,C_YZ,U,W)\\
&-B(Y,C_XZ,U,W),
\end{split}
\end{equation}
and
\begin{equation}\label{app:vector-AC}
\begin{split}
0={}&B(Y,C_XZ,U,W)-B(X,C_YZ,U,W)\\
&-B(Y,C_XU,Z,W)+B(X,C_YU,Z,W).
\end{split}
\end{equation}
The five arguments here are arbitrary tangent vectors. These formulas distinguish vector arguments from the component indices used in the canonical tensor calculation at the end of the appendix.

\begin{definition}[Induced centroaffine connection]
For tangent vector fields $X,Y$ on $\Sigma$, define $\nablac$ by
\[
D_XY=\nablac_XY-h(X,Y)y.
\]
Thus the transversal vector field is $-y$, and the induced tangent connection is $\nablac=\nablah-C^\sharp$.
\end{definition}

\begin{lemma}\label{geom:identities}
The connection $\nablac$ is torsion-free and satisfies
\[
\nablac h=2C,\qquad
R^{\nablac}(X,Y)Z=h(Y,Z)X-h(X,Z)Y.
\]
The tensor $\nablah C$ is fully symmetric.
\end{lemma}
\begin{proof}
The connection formula follows from~\eqref{eq:angularD}. Since $C^\sharp(X,Y)=C^\sharp(Y,X)$, $\nablac$ is torsion-free. Substitution into $\nablah h=0$ gives $\nablac h=2C$. Taking the tangent component of $R^D(X,Y)Z=0$ gives the displayed curvature identity. Full symmetry of $\nablah C$ was established in Section~\ref{sec:preliminaries}.
\end{proof}

\begin{lemma}[Closed centroaffine geodesics]\label{geom:periodic}
For $y_0\in\Sigma$ and $0\ne w_0\in T_{y_0}\Sigma$, the $\nablac$-geodesic with those initial data is defined for all real values of the affine parameter, lies in $\Sigma\cap\operatorname{span}(y_0,w_0)$, and has periodic position and velocity. Its speed squared $a(t)=h(w(t),w(t))$ is positive, periodic, and bounded above and away from zero.
\end{lemma}
\begin{proof}
Orient the plane $\Pi=\operatorname{span}(y_0,w_0)$ so that $c=\det(y_0,w_0)>0$, using any auxiliary Euclidean coordinates. Parametrize the radial section by
\[
\eta(\theta)=r(\theta)(\cos\theta,\sin\theta),\qquad
r(\theta)=F(\cos\theta,\sin\theta)^{-1}.
\]
The function $r$ is smooth, positive, and $2\pi$-periodic, and $\det(\eta,\eta_\theta)=r^2$. Define affine time by
\[
\frac{\dd t}{\dd\theta}=\frac{r(\theta)^2}{c},
\]
with its additive constant chosen at $y_0$. This is a global increasing parametrization of the real line, with a positive finite period increment. Its velocity has the prescribed initial value because both tangent vectors have determinant $c$ with $y_0$.

The determinant $\det(y,\dot y)=c$ is constant, so $\ddot y$ is radial. Differentiate $F(y)=1$ twice. On angular vectors $D^2F=h$, and $\dd F_y(y)=1$, giving
\[
\ddot y=-h(\dot y,\dot y)y.
\]
The preceding lemma identifies this with $\nablac_{\dot y}\dot y=0$. Smoothness, periodicity and nonvanishing of the velocity prove all remaining assertions.
\end{proof}

Suppose first that
\[
0<q:=\max_{\Sigma}\rank\{v\mapsto B(v,\cdot,\cdot,\cdot)\}<m,
\qquad
U:=\{y:\rank B_y=q\}.
\]
The rank here is the rank of the insertion map, not the rank of the quartic
viewed as an operator on symmetric two-tensors. The set $U$ is nonempty and
open: an order-$q$ nonzero minor remains nonzero locally, while maximality
prevents a larger rank.

Write $\nablah$ for the Levi--Civita connection of $h$, and let
$\nablac=\nablah-C^\sharp$ be the induced centroaffine connection.
We use the identities established earlier,
\[
\nablac h=2C,\qquad
R^{\nablac}(X,Y)Z=h(Y,Z)X-h(X,Z)Y,\qquad
\nablah C\in\Sym^4T^*\Sigma.
\]
By Lemma~\ref{geom:periodic}, every nonconstant $\nablac$-geodesic is a complete periodic parametrization of a planar section.

On $U$ let
\[
N=\operatorname{rad}B,\qquad P=N^{\perp_h},\qquad
b=B|_{P^4}.
\]
Constant rank makes $N$ and $P$ smooth bundles of ranks $m-q$ and $q$. The orthogonal complement $P$ is called the active subbundle of $B$ on $U$. Every value of $B$ is determined by projecting its arguments
onto $P$, and $b$ has trivial radical. Indeed, if $X\in P$ annihilates
$b(X,\cdot,\cdot,\cdot)$, then the projection identity makes it annihilate
$B$ on the full tangent space; thus $X\in N\cap P=0$.

\begin{definition}[Centroid of a symmetric quartic]
For a quartic form $b$ with trivial radical on $P$, its centroid is
\[
\Cent(b)=\{A\in\End(P):
 b(AX,Y,Z,W)=b(X,AY,Z,W)\ \text{for all arguments}\}.
\]
Thus an element of the centroid can be transferred between any two argument slots by symmetry.
\end{definition}

\begin{lemma}\label{deg:centroid}
The space $\operatorname{Cent}(b)$ is a commutative unital algebra under
composition.
\end{lemma}
\begin{proof}
The defining identity is linear in the endomorphism and is satisfied by $\Id$. If $A,D$ belong to it,
successive applications of the defining identity give
\begin{align*}
b(ADX,Y,Z,W)
 &=b(DX,AY,Z,W)=b(X,AY,DZ,W)\\
 &=b(AX,Y,DZ,W)=b(DAX,Y,Z,W).
\end{align*}
Triviality of the radical gives $AD=DA$. A further application of the
defining identity gives
\[
b(ADX,Y,Z,W)=b(DAX,Y,Z,W)
=b(AX,DY,Z,W)=b(X,ADY,Z,W).
\]
Thus $AD$ also belongs to the centroid.
\end{proof}

We first extract two consequences of the compatibility identity
\eqref{app:vector-AC}. Insert $X=u\in N$ and $U=v\in N$, with $Y,Z,W\in P$.
All terms vanish except the insertion of $C^\sharp(u,v)$ into $b$.
Its active projection is therefore zero. Hence
\begin{equation}\label{deg:mixed}
C(N,N,P)=0,\qquad C^\sharp(N,N)\subset N,\qquad C^\sharp(N,P)\subset P.
\end{equation}
For $v\in N$ define
\[
A_vX=C^\sharp(v,X)\qquad(X\in P).
\]
Inserting $X=v$ and taking the other arguments in $P$ in
\eqref{app:vector-AC} shows that $A_v$ satisfies the centroid condition. Hence
\begin{equation}\label{deg:Acent}
A_v\in\operatorname{Cent}(b),\qquad A_v=A_v^*.
\end{equation}
Self-adjointness follows separately from
$h(A_vX,Y)=C(v,X,Y)=h(X,A_vY)$.

Next take $X=u$, $Z=v$ in $N$ in the twisted Codazzi equation
\eqref{app:vector-TC}, with the other arguments in $P$. Its right side is zero
by~\eqref{deg:mixed}. Terms whose derivatives contain two radical arguments vanish.
Expanding the derivative with one radical slot leaves
\[
b\bigl(Y,(\nablah_uv)_P,U,W\bigr)=0.
\]
Thus
\begin{equation}\label{deg:autoparallel}
\nablah_NN\subset N,\qquad \nablah_NP\subset P.
\end{equation}
The first inclusion and torsion-freeness make $N$ integrable; the second
follows by differentiating orthogonality. Combining
\eqref{deg:autoparallel} with~\eqref{deg:mixed} gives the corresponding
inclusions for $\nablac$:
\begin{equation}\label{deg:affine-preservation}
\nablac_NN\subset N,\qquad \nablac_NP\subset P.
\end{equation}

For $v\in N$ and $X\in P$ put
\begin{equation}\label{deg:splittings}
T_vX=-(\nablah_Xv)_P,\qquad
S_v=T_v+A_v=-(\nablac_Xv)_P.
\end{equation}
These are tensorial in the pointwise values of $v$ and $X$: differentiating
a scalar multiplying $v$ contributes only a radical vector and is removed
by the projection.

Take the differentiation direction $v\in N$ and all four quartic arguments in $P$ in~\eqref{app:vector-TC}.
Here $\nabla^{h,P}$ denotes the connection induced by $\nablah$ on $P$
along radical directions; preservation in \eqref{deg:autoparallel}
makes it well defined.
The derivative of the identically zero insertion $B(v,Y,Z,W)$ is
$-B(\nablah_Xv,Y,Z,W)=b(T_vX,Y,Z,W)$. Consequently,
\begin{equation}\label{deg:quartic-h}
(\nabla^{h,P}_vb)(X,Y,Z,W)=b((T_v-A_v)X,Y,Z,W).
\end{equation}
Its left side is symmetric in the four arguments. Thus $T_v-A_v$
belongs to the centroid. By~\eqref{deg:Acent} and
Lemma~\ref{deg:centroid}, all operators
\[
A_v,\ T_v,\ S_v\qquad(v\in N)
\]
belong to the same commutative algebra.

Choose $y_0\in U$ and $0\ne v_0\in N_{y_0}$. Let $y(t)$ be the complete
$\nablac$-geodesic with these initial data, and set
\[
v(t)=y'(t),\qquad a(t)=h(v(t),v(t)).
\]
Lemma~\ref{geom:periodic} gives periodic position and velocity, and constants
$0<a_0\leq a(t)\leq a_1$ on the real line.

While the curve stays in $U$, its velocity stays radical. The assertion is
locally true by~\eqref{deg:affine-preservation} and geodesic uniqueness.
It is also a closed condition, since $B(v,\cdot,\cdot,\cdot)=0$ is preserved
under limits. An open-and-closed argument on the residence interval
therefore proves the assertion throughout that interval.

Along this radical curve, take a $\nablac$-parallel frame of $P$. By the
ambient centroaffine formula,
\[
D_vX=\nablac_vX-h(v,X)y=0
\]
for every such parallel active vector. Therefore this is a frame of
\emph{fixed vectors in the original space $V$}. In particular,
$P_{y(t)}=P_0:=P_{y_0}$ throughout the residence interval.

\begin{lemma}\label{deg:matrix-system}
In this fixed frame, the splitting matrix $S=S_{v(t)}$ and the raised
quartic insertion map $M$ satisfy
\begin{equation}\label{deg:ODE}
S'=S^2+a(t)\Id,\qquad M'=MS,
\end{equation}
and $M$ is injective.
\end{lemma}
\begin{proof}
The quotient description of the splitting operator is
\[
S_v[X]=-[\nablac_Xv]\quad\text{on }T\Sigma/N.
\]
It is independent of the representative by~\eqref{deg:affine-preservation}.
Extend $v$ locally to a radical geodesic field. One constructs this extension
by taking geodesics with radical initial velocities through a transversal;
integrability and the autoparallel property keep the short geodesics in their
radical leaves. Choose a local field $X$ with $[v,X]=0$. Torsion-freeness
then gives $\nablac_v[X]=-S_v[X]$. The quotient curvature calculation is
\[
[R^{\nablac}(v,X)v]=-(\nablac_vS_v)[X]+S_v^2[X].
\]
On the other hand the centroaffine curvature identity gives
$[R^{\nablac}(v,X)v]=-a(t)[X]$. Thus
$\nablac_vS_v=S_v^2+a\Id$, which is the first equation of~\eqref{deg:ODE}
in the parallel active frame.

Define the raised tensor $b^\sharp$ by
\[
h(b^\sharp(X,Y,Z),W)=b(X,Y,Z,W),
\qquad
(MX)(Y,Z)=b^\sharp(X,Y,Z).
\]
Its insertion map $M:P_0\to\operatorname{Hom}(\Sym^2P_0,P_0)$ is injective
because $b$ has trivial radical.

Write $\nabla^{\mathrm c,P}$ for the connection on $P$ induced by
$\nablac$ along radical directions; it is well defined by
\eqref{deg:affine-preservation}. Changing the connection in~\eqref{deg:quartic-h} adds four $A_v$ terms,
one in each covariant slot. Their centroid property yields
\[
(\nabla^{\mathrm c,P}_vb)(X,Y,Z,W)
=b((T_v+3A_v)X,Y,Z,W)
=b((S_v+2A_v)X,Y,Z,W).
\]
Meanwhile $\nablac_vh=2h(A_v\cdot,\cdot)$. Self-adjointness and the centroid
property imply
\[
b^\sharp(A_vX,Y,Z)=A_vb^\sharp(X,Y,Z).
\]
Differentiating the defining metric contraction therefore cancels the
two $A_v$ terms:
\[
\nabla^{\mathrm c,P}_vb^\sharp(X,Y,Z)=b^\sharp(S_vX,Y,Z).
\]
With parallel arguments in the fixed linear frame this is $M'=MS$.
\end{proof}

Let $f,g$ solve
\[
z''+a(t)z=0,\qquad
(f(0),f'(0))=(1,0),\qquad(g(0),g'(0))=(0,1),
\]
and put $Y=f\Id-gS_0$. Their Wronskian is $fg'-f'g=1$. In particular
$f,g$ never vanish simultaneously. Direct differentiation of the inverse
matrix shows that, while $Y$ is invertible,
\begin{equation}\label{deg:denominator}
S=-Y'Y^{-1},\qquad M=M_0Y^{-1}
\end{equation}
solve~\eqref{deg:ODE} with the prescribed initial data. The uniqueness theorem
for ordinary differential equations identifies these expressions with the
geometrically defined operators on their common interval.

\begin{lemma}[Rank retention and spectral obstruction]\label{deg:no-real}
The curve cannot leave $U$ at a finite time before the first singularity
of $Y$. Moreover, $S_v$ has no real eigenvalue for any $0\ne v\in N$.
\end{lemma}
\begin{proof}
Suppose the curve has a finite residence endpoint $t_0$ with $Y(t_0)$
invertible. Its fixed ambient active frame remains tangent at the endpoint,
by continuity of $\dd E_{y(t)}(X)=0$. Restrict the global smooth tensor $B$
and the positive definite endpoint metric to this fixed space $P_0$.
Raising an index there gives exactly the finite limit of $M(t)$.
Equation~\eqref{deg:denominator} identifies that limit as
$M_0Y(t_0)^{-1}$, which is injective. Thus the restricted endpoint quartic
has insertion rank $q$, so the full endpoint rank is at least $q$.
Global maximality makes it exactly $q$; hence the endpoint lies in $U$.
Continuity also keeps the velocity radical there. Local existence then
continues the residence interval, a contradiction.

Suppose $S_0$ has a real eigenvalue $\lambda$. The scalar
$u=f-\lambda g$ starts at $u(0)=1$ and satisfies $u''+au=0$. It has a finite
positive zero. Indeed, if it stayed positive, then $w=u'/u$ would satisfy
\[
w'=-a-w^2\leq-a_0-w^2,
\qquad
\frac{\dd}{\dd t}\arctan\!\frac{w}{\sqrt{a_0}}
\leq-\sqrt{a_0},
\]
which is impossible for all positive time.

Thus $Y$ has a first positive singular time $t_*$. At that time $g(t_*)$
cannot be zero, and $\mu=f(t_*)/g(t_*)$ is a real eigenvalue of $S_0$.
Choose a nonzero real eigenvector $Z$. For $0\leq t<t_*$,
rank retention ensures that the geometrically defined operators satisfy~\eqref{deg:denominator}, so
\[
M(t)Z=\frac{M_0Z}{f(t)-\mu g(t)}.
\]
The numerator is nonzero by injectivity and the denominator tends to zero.
This diverges, whereas restriction to the fixed subspace $P_0$ gives a
finite limit for $M(t)$ at $t_*$, independently of the endpoint rank. This contradiction excludes every real eigenvalue.
\end{proof}

For later use, the same argument also proves that the whole radical
geodesic stays in $U$ for all real times. Indeed the absence of real
eigenvalues makes $Y(t)$ invertible for all real $t$: if $g(t)=0$, use
$f(t)\ne0$; otherwise a singularity would give the real eigenvalue
$f(t)/g(t)$. Rank retention excludes every finite residence endpoint,
in either time direction.

At a fixed point the family $\{S_v:v\in N_y\}$ is a commuting family of
real operators, real-linear in $v$. After complexification it has a common
eigenvector $Z\ne0$. Indeed, this follows by successively choosing
eigenspaces of a finite set of commuting generators; each chosen nonzero
eigenspace is invariant under the remaining generators. Write
\[
S_vZ=\ell(v)Z,\qquad\ell:N_y\longrightarrow\mathbb C
\]
for the resulting real-linear eigenvalue functional.
If $\dim N_y\geq2$, a nonzero $v$ lies in $\ker\operatorname{Im}\ell$.
Then $S_v$ has a real eigenvalue; taking a nonzero real or imaginary part
of $Z$ gives a real eigenvector. This contradicts Lemma~\ref{deg:no-real}.
Thus $\dim N=1$ and $q=m-1$. A real matrix of odd order has a real eigenvalue,
so $q$ is even.

Choose a local $h$-unit section $v$ of the radical line and write
\[
c=C(v,v,v),\qquad A=A_v,\quad T=T_v,\quad S=T+A.
\]
Unit length and~\eqref{deg:autoparallel} imply $\nablah_vv=0$.
Equation~\eqref{deg:mixed} gives $C^\sharp(v,v)=cv$.
For $X\in P$, unit length also gives $\nablah_Xv=-TX$ without any
radical correction.

In an $h$-parallel active frame along a unit radical curve, differentiation
of $h(AX,Y)=C(v,X,Y)$ and the Cartan--Codazzi equation give
\begin{align*}
h(A'X,Y)
 &=(\nablah_vC)(v,X,Y)
   =(\nablah_XC)(v,v,Y)\\
 &=2h(ATX,Y)-c\,h(TX,Y).
\end{align*}
Here the first two slot corrections produce $2A$, while the last uses
$h(v,\nablah_XY)=h(TX,Y)$. Therefore
\begin{equation}\label{deg:A-derivative}
A'=(2A-c\Id)T.
\end{equation}
Both $A$ and $A'$ are self-adjoint. Also $AT=TA$ by the centroid algebra.
On an eigenspace $E_\alpha$ of $A$, if $2\alpha-c\ne0$,
self-adjointness of~\eqref{deg:A-derivative} makes $T|_{E_\alpha}$
self-adjoint. Then $S=T+A$ has a real eigenvalue, contradicting
Lemma~\ref{deg:no-real}. Consequently
\[
A=(c/2)\Id,\qquad A'=0,\qquad v(c)=0.
\]

Take the complete periodic radical centroaffine geodesic through any
point of $U$. Write $w=y'=uv$ with $u=|w|_h>0$ and $v=w/u$.
The preceding argument applies consistently along this orbit, without an
orientation of $N$ on all of $U$, and shows that $c$ is constant along it.
The two connections give
\[
\nablah_ww=C^\sharp(w,w)=u^2cv,
\qquad
\nablah_ww=u'v,
\]
so $u'=cu^2$. With respect to $h$-arc length $s$, where $\dd s/\dd t=u$, this becomes
$\dd(\log u)/\dd s=c$. Since the velocity is periodic, $u$ returns to its initial value after one period.
Integration over its positive $h$-length gives $c=0$. We have proved
\begin{equation}\label{deg:radical-C-zero}
C(N,T\Sigma,T\Sigma)=0,\qquad A_v=0,\qquad S_v=T_v
\quad\text{on }U.
\end{equation}

Let $\gamma=C|_{P^3}$. Using~\eqref{deg:radical-C-zero} in the
Cartan--Codazzi equation gives, for active arguments,
\[
(\nabla^{h,P}_v\gamma)(X,Y,Z)
=(\nablah_XC)(v,Y,Z)=\gamma(T_vX,Y,Z).
\]
The left side is symmetric in $X,Y,Z$, so $T_v$ can be moved between any
two slots of $\gamma$. It already has that property for $b$ and has no
real spectrum.

Fix one point and a nonzero radical vector. Complexify the active space
and $T=T_v$. Let $P_+$ and $P_-$ be the sums of its generalized eigenspaces
with eigenvalues respectively in the upper and lower half-planes. Then
\[
P_{\mathbb C}=P_+\oplus P_-,
\qquad \overline{P_+}=P_-.
\]
The primary factors of the minimal polynomial for these two halves are
coprime. Bezout's identity constructs their projections as complex
polynomials in $T$. Every such polynomial can be moved between slots of
$b$ and $\gamma$. Moving the projection onto $P_+$ from an argument in $P_+$
to an argument in $P_-$ proves that all mixed components vanish. Thus both tensors have only all-positive and all-negative parts.
The restriction $b_-$ to $P_-^4$ has trivial radical: a vector annihilating
it also annihilates the full complexified $b$, whose insertion map is the
complexification of an injective real map.

Apply~\eqref{app:vector-AC} with its first and third arguments in $P_+$ and the remaining three in $P_-$.
The second and third terms vanish because the mixed cubic components vanish;
the fourth term has mixed fixed arguments in $b$. Hence
\[
b(w_1,\gamma^\sharp(v_1,v_2),w_2,w_3)=0
\quad(v_1,v_2\in P_+,\ w_1,w_2,w_3\in P_-).
\]
Set $\alpha=\gamma(v_1,v_2,\cdot)$ and $Z=h_{\mathbb C}^{-1}\alpha$.
Triviality of the radical of $b_-$ makes the $P_-$ component of $Z$ zero.
Thus $Z\in P_+$. On the other hand $\alpha$ annihilates $P_-$, so its
conjugate covector $\bar\alpha$ annihilates $P_+$. Positive definiteness gives
\[
0=\bar\alpha(Z)
=h_{\mathbb C}^{-1}(\alpha,\bar\alpha)
=|\operatorname{Re}\alpha|_{h^{-1}}^2+
 |\operatorname{Im}\alpha|_{h^{-1}}^2.
\]
Therefore $\alpha=0$. The positive part of $\gamma$ vanishes, and conjugation
eliminates its negative part. Thus $\gamma=0$.

\begin{proposition}[The deficient-rank case]
\label{degenerate:conclusion}
If $0<q<m$, then $q=m-1$ is even and $C=0$ on the entire maximal-rank set $U$.
\end{proposition}
\begin{proof}
The preceding argument gives $\dim N=1$ and shows that $q=m-1$ is even. The vanishing of the transverse cubic,
together with~\eqref{deg:radical-C-zero}, eliminates every component of $C$.
\end{proof}

\begin{definition}[Irreducible family of endomorphisms]
A family of endomorphisms of a finite-dimensional vector space $H$ is irreducible if its only common invariant linear subspaces are $\{0\}$ and $H$. Invariance of $L\subset H$ means $A(L)\subset L$ for every member $A$ of the family.
\end{definition}

Put
\[
\Omega=\{p\in\Sigma:\operatorname{rad}B_p=0\}.
\]
For the full-rank case, assume that $\Omega\ne\varnothing$. For $u,v\in T_p\Sigma$
define the self-adjoint contraction $B_{uv}$ by
\[
h(B_{uv}z,w)=B(u,v,z,w).
\]
The self-adjoint common commutant is
\[
\mathcal Z_B(p)=
\{S\in\End(T_p\Sigma):S=S^*,\ [S,B_{uv}]=0\text{ for all }u,v\in T_p\Sigma\}.
\]
It contains $\R\Id$. We will prove that equality holds on a nonempty open
dense subset of $\Omega$. The proof is local.

\begin{lemma}\label{fullrank:no-splitting}
There is no nonempty open set $O\subset\Omega$ with smooth orthogonal
subbundles $P,Q$ of positive rank such that
\[
T\Sigma|_O=P\oplus Q,\qquad B=B_P+B_Q.
\]
Here $B_P$ and $B_Q$ are the restrictions to $P^4$ and $Q^4$, extended by zero
on the other summand. Thus the asserted equality means that every mixed
component of $B$ is zero.
\end{lemma}
\begin{proof}
Suppose that such a splitting exists. Each restricted quartic has trivial
radical. Indeed, if $w\in P_p$ satisfies $B_P(w,P_p,P_p,P_p)=0$, expansion
of the last three arguments according to $P_p\oplus Q_p$ gives
$B_p(w,\cdot,\cdot,\cdot)=0$. Since $p\in\Omega$, $w=0$. The argument
for $B_Q$ is identical.

We first show that the cubic also splits. In the compatibility identity
\eqref{app:vector-AC}, take the first and fourth arguments in $Q$ and the remaining three in $P$. All terms except
the one containing the Cartan contraction of the two $Q$ vectors have fixed arguments from both
summands and vanish. Therefore
\[
B_P\bigl(X,(C^\sharp(Y,Z))_P,U,W\bigr)=0
\quad(Y,Z\in Q,\ X,U,W\in P).
\]
The values of $X,U,W$ are arbitrary. Triviality of the radical of $B_P$
implies $(C^\sharp(Y,Z))_P=0$, or $C(Q,Q,P)=0$. Interchanging $P$ and $Q$
gives $C(P,P,Q)=0$. Symmetry of the cubic now yields
\begin{equation}\label{fullrank:C-split}
C=C_P+C_Q.
\end{equation}
Consequently, for $X\in P$, $C_X$ maps $P$ to itself and is zero on $Q$;
for $Y\in Q$, $C_Y$ maps $Q$ to itself and is zero on $P$. In particular,
\begin{equation}\label{fullrank:mixed-commutator}
[C_X,C_Y]=0\qquad(X\in P,\ Y\in Q).
\end{equation}

We next show that both subbundles are parallel for $\nabla^h$. Choose smooth local sections $X_1,X_2,X_3$ of $P$ and $Y,Z$
of $Q$. Apply~\eqref{app:vector-TC} with differentiation directions $Y,X_1$
and the remaining arguments $X_2,X_3,Z$. This gives
\begin{equation}\label{fullrank:mixed-TC}
(\nablah_YB)(X_1,X_2,X_3,Z)
=(\nablah_{X_1}B)(Y,X_2,X_3,Z).
\end{equation}
Both terms on its original right-hand side are zero: irrespective of the
contracted vector, their other arguments include vectors from both blocks.

The expression on the right of \eqref{fullrank:mixed-TC} is zero. Its scalar
derivative differentiates an identically zero mixed component. Every connection
correction replaces just one argument of a component with two arguments in each
block; the three unchanged arguments still include both blocks. On the left,
the same expansion leaves exactly one potentially nonzero correction:
\[
(\nablah_YB)(X_1,X_2,X_3,Z)
=-B_P\bigl(X_1,X_2,X_3,(\nablah_YZ)_P\bigr).
\]
Indeed, correcting one of the $X_r$ leaves two unchanged $P$ arguments and the
unchanged $Q$ argument $Z$. It follows from triviality of the radical of $B_P$
that
\[
\nablah_QQ\subset Q.
\]
Interchanging $P$ and $Q$ gives $\nablah_PP\subset P$. Orthogonality and metric
compatibility give the cross-direction inclusions as well. For example,
for $X\in P$ and $Y,Z\in Q$,
\[
h(\nablah_YX,Z)=-h(X,\nablah_YZ)=0,
\]
so $\nablah_QP\subset P$. Similarly $\nablah_PQ\subset Q$. Thus both
subbundles are preserved by $\nablah$ in every direction.

At any point of $O$ choose orthonormal vectors $X\in P$ and $Y\in Q$.
Since $Q$ is parallel, $R^h(X,Y)Y\in Q$, and hence
\[
h(R^h(X,Y)Y,X)=0.
\]
On the other hand, the angular Gauss identity \eqref{eq:Gauss} and
\eqref{fullrank:mixed-commutator} give
\[
R^h(X,Y)Y=R_0(X,Y)Y=X.
\]
These identities therefore yield $0=1$, a contradiction.
\end{proof}

\begin{remark}
Smoothness of the splitting on an open set is essential. A pointwise
orthogonal decomposition cannot be differentiated. The next argument produces
the needed smooth subbundles from a nontrivial self-adjoint common commutant
on an open set.
\end{remark}

\begin{proposition}\label{fullrank:generic}
The set
\[
\mathcal G=\{p\in\Omega:\mathcal Z_B(p)=\R\Id\}
\]
is open and dense in $\Omega$. In particular, it is nonempty.
\end{proposition}
\begin{proof}
On a coordinate neighborhood in $\Omega$ choose a smooth $h$-orthonormal frame
$e_1,\ldots,e_m$. The bundle $\mathcal E$ of self-adjoint endomorphisms has rank
$d=m(m+1)/2$. Consider the smooth linear bundle map
\[
\mathcal L_p:\mathcal E_p\longrightarrow
\bigoplus_{i\le j}\mathfrak{so}(T_p\Sigma,h_p),
\qquad
S\longmapsto([S,B_{e_ie_j}])_{i\le j}.
\]
Its kernel is $\mathcal Z_B(p)$, because the finitely many $B_{e_ie_j}$
span all contractions by bilinearity. Since $\Id$ belongs to the kernel,
$\rank\mathcal L_p\le d-1$.

If $p\in\mathcal G$, some minor of order $d-1$ is nonzero at $p$ and
remains nonzero nearby. The rank cannot exceed $d-1$, so it equals $d-1$
on that neighborhood. This proves openness.

For density, suppose that some nonempty coordinate open set
$O\subset\Omega$ contains no point of $\mathcal G$. Then
$\rank\mathcal L\le d-2$ throughout $O$. Choose a point where this rank is
maximal on $O$, and shrink to a neighborhood $U$ on which it is constant.
This is possible because rank takes finitely many integer values and a
nonzero maximal minor persists. If the maximal rank is zero, any sufficiently
small neighborhood already has constant rank.

The kernels on $U$ form a smooth vector bundle of rank at least two.
A fixed invertible maximal minor determines the pivot variables smoothly in
terms of the free variables, yielding a smooth local kernel frame. Choose a nonscalar element in one fibre and
extend it to a smooth kernel section $S$. Replacing it by
$S-(\tr S/m)\Id$ and shrinking $U$ gives a nowhere-zero traceless section.
It is therefore nonscalar at every point, self-adjoint, and commutes with
every $B_{uv}$.

At one point $p_0$, select a nonempty proper cluster of eigenvalues of $S(p_0)$,
separated from the others by a positive gap. Such a cluster exists because
$S(p_0)$ is self-adjoint and nonscalar. The separation persists after shrinking $U$ again. A fixed closed contour $\gamma$ in the complex plane can
be chosen to surround the selected cluster and no other eigenvalues at any
point of $U$. Set
\[
\Pi(p)=\frac{1}{2\pi\mathrm i}
\int_\gamma(z\Id-S(p))^{-1}\,\dd z.
\]
Smoothness follows from smoothness of the inverse matrix away from its
singular set and differentiation under this fixed contour integral.
Diagonalization of $S(p)$ shows that $\Pi(p)$ is the real orthogonal
projection onto the selected eigenspaces. Its rank
is locally constant, positive, and strictly less than $m$.

Commutation with $S$ implies commutation with its resolvent, hence with
$\Pi$. Consequently $P=\im\Pi$ and $Q=\ker\Pi$ are smooth positive-rank
orthogonal subbundles invariant under every contraction $B_{uv}$. For
$z\in P$ and $w\in Q$,
\[
B(u,v,z,w)=h(B_{uv}z,w)=0.
\]
Complete symmetry moves any two arguments belonging to different blocks into
these last two positions. All mixed components are thus zero:
$B=B_P+B_Q$. This contradicts Lemma~\ref{fullrank:no-splitting}, proving
density.
\end{proof}

The equality $\mathcal Z_B(p)=\R\Id$ is equivalent to irreducibility of the contraction family. A proper common
invariant subspace has an invariant orthogonal complement, since all contractions
are self-adjoint; its orthogonal projection is a nonscalar element of the common commutant.
Conversely, a nonscalar self-adjoint element of the common commutant has a proper eigenspace
preserved by every contraction. Thus Proposition~\ref{fullrank:generic} also
proves generic irreducibility.

For the pointwise algebraic argument below, let $H$ be a finite-dimensional real inner-product space, of
dimension $m\geq2$, with inner product $h$. Let
\[
 C\in\operatorname{Sym}^3H^*,\qquad
 B\in\operatorname{Sym}^4H^*,\qquad B\neq0.
\]
Define self-adjoint endomorphisms by
\[
 h(C_xy,z)=C(x,y,z),\qquad
 h(B_{uv}x,y)=B(u,v,x,y),
\]
and put $K(x,y)=[C_x,C_y]$. The angular compatibility identity is
\begin{equation}\label{alg:F}
 \begin{split}
 0={}&B(u,C_xv,w,z)-B(x,C_uv,w,z)\\
    &-B(u,C_xw,v,z)+B(x,C_uw,v,z)
 \end{split}
\end{equation}
for all $x,u,v,w,z\in H$. This is the pointwise form of~\eqref{eq:angularCompatibility}.

\begin{theorem}\label{algebra:isotropy}
Suppose \eqref{alg:F} holds and
\[
 \{A=A^*:[A,B_{uv}]=0\text{ for every }u,v\}
   =\R\Id.
\]
Then there is a real scalar $s$ such that
\[
 [C_x,C_y]z=s\{h(y,z)x-h(x,z)y\}.
\]
\end{theorem}

We prove the theorem in the following lemmas, beginning with a consequence of~\eqref{alg:F}.

The tensor $K$ has the algebraic curvature symmetries. Indeed, $K(x,y)$ is
skew-adjoint, and
\[
 K(x,y)z+K(y,z)x+K(z,x)y=0
\]
follows by expanding and using $C_xy=C_yx$. Pair interchange also follows
directly from
\[
 h(K(x,y)z,w)
 =h(C_yz,C_xw)-h(C_xz,C_yw).
\]
We shall use the four-tensor convention
$R(x,y,z,w)=h(K(x,y)z,w)$.

\begin{lemma}\label{alg:quartic-compatible}
Every contraction $B_{uv}$ satisfies
\begin{equation}\label{alg:K-compatible}
 K(x,y)B_{uv}z+K(y,z)B_{uv}x+K(z,x)B_{uv}y=0.
\end{equation}
\end{lemma}
\begin{proof}
Set
\[
 T(x,u,v)=[C_x,B_{uv}].
\]
Pairing the difference $T(x,u,z)-T(u,x,z)$ with arbitrary $v,w$ shows that
\[
 h\bigl((T(x,u,z)-T(u,x,z))v,w\bigr)
 =-\mathcal F(x,u,v,w,z)=0,
\]
where $\mathcal F$ is the four-term expression in \eqref{alg:F}.
Thus \eqref{alg:F} says $T(x,u,z)=T(u,x,z)$. Since $B_{uv}=B_{vu}$,
$T$ is symmetric in all three of its displayed arguments.

For fixed $u,v$ write
\[
 S_{uv}(x,y,z)=\sum_{\cyc(x,y,z)}[C_x,C_y]B_{uv}z.
\]
Commuting the innermost $C$ past $B_{uv}$, and then using $C_yz=C_zy$,
gives
\begin{align}
 S_{uv}
 &=\sum_{\cyc}C_x\{T(y,u,v)z-T(z,u,v)y\}.
 \label{alg:S-first}
\end{align}
Here the discarded terms are
$\sum_{\cyc}(C_xB_{uv}C_yz-C_yB_{uv}C_xz)$, which cancel pairwise.
The full symmetry of $T$ and that of $B$ imply
\begin{align*}
 T(y,u,v)z-T(z,u,v)y
 &=[C_u,B_{yv}]z-[C_u,B_{zv}]y\\
 &=B_{zv}C_uy-B_{yv}C_uz,
\end{align*}
because $B_{yv}z=B_{zv}y$.
Substitute this into \eqref{alg:S-first} and commute $C_x$ past the two
contractions. The extra terms
\[
 \sum_{\cyc}
 \{T(x,z,v)C_uy-T(x,y,v)C_uz\}
\]
again cancel pairwise. Regrouping the remaining cyclic terms yields
\begin{equation}\label{alg:S-switch}
 S_{uv}(x,y,z)
 =\sum_{\cyc}B_{zv}(C_xC_uy-C_yC_ux)
 =\sum_{\cyc}B_{zv}K(x,y)u.
\end{equation}
For the last equality, use $C_uy=C_yu$ and $C_ux=C_xu$.

Fix $x,y,z$ temporarily and put $t(u,v,w)=h(S_{uv}(x,y,z),w)$.
Equation \eqref{alg:S-switch}, full symmetry of $B$, and skew-adjointness of
$K(x,y)$ give
\begin{align*}
 t(u,v,w)
 &=\sum_{\cyc}B(z,v,K(x,y)u,w)\\
 &=-h\left(u,\sum_{\cyc}K(x,y)B_{vw}z\right)
 =-t(v,w,u).
\end{align*}
Applying this relation three times gives
$t(u,v,w)=-t(v,w,u)=t(w,u,v)=-t(u,v,w)$.
Thus $t=0$ over $\R$. Since $w$ is arbitrary, $S_{uv}=0$.
This proves \eqref{alg:K-compatible}.
\end{proof}

\begin{definition}[Curvature-compatible self-adjoint operator]
A self-adjoint operator $A$ is compatible with the algebraic curvature tensor $K$ if
\[
 \mathcal C_A(x,y,z)
 :=K(x,y)Az+K(y,z)Ax+K(z,x)Ay=0
\]
for all $x,y,z$. Compatibility is linear in $A$, and the identity operator is compatible by the first Bianchi identity.
\end{definition}

\begin{lemma}\label{alg:jordan-closure}
If $A,D$ are compatible self-adjoint operators, then
$A\circ D=(AD+DA)/2$ is compatible.
\end{lemma}
\begin{proof}
Define
\[
 U(x,y,z,w)
 =R(x,y,Az,Dw)+R(x,y,Dz,Aw).
\]
It is skew in the first pair and in the second pair.
Its cyclic sum over $x,y,z$ vanishes, because the two sums are respectively
$h(\mathcal C_A(x,y,z),Dw)$ and
$h(\mathcal C_D(x,y,z),Aw)$.
Hence $U$ has the algebraic curvature symmetries. Pair interchange follows from skewness in each pair and the first Bianchi identity by the same four-tensor calculation as the one following~\eqref{round:algebraic-bianchi}.

Expand the compatibility expression for $AD+DA$. For example,
compatibility of $A$ at $(x,y,Dz)$ gives
\[
 R(x,y,ADz,w)
 =-R(y,Dz,Ax,w)-R(Dz,x,Ay,w).
\]
Use the analogous formula with $A,D$ interchanged and take the cyclic sum.
After regrouping the terms according to their first argument, the result is
\begin{align*}
 h(\mathcal C_{AD+DA}(x,y,z),w)
 =\sum_{\cyc}\{&
 -R(z,Dx,Ay,w)+R(z,Ay,Dx,w)\\
 &-R(z,Ax,Dy,w)+R(z,Dy,Ax,w)\}.
\end{align*}
The Bianchi identity for $R$ turns each brace into
\[
 R(Dx,Ay,z,w)+R(Ax,Dy,z,w).
\]
Pair interchange for $R$ therefore gives
\[
 h(\mathcal C_{AD+DA}(x,y,z),w)
   =\sum_{\cyc}U(z,w,x,y).
\]
Pair interchange for $U$ turns this into
$\sum_{\cyc}U(x,y,z,w)=0$.
Thus $AD+DA$, and hence $A\circ D$, is compatible.
\end{proof}

\begin{definition}[Jordan algebra generated by the quartic contractions]
For self-adjoint operators define $A\circ D=(AD+DA)/2$. A real linear subspace closed under this product is a Jordan subalgebra of the self-adjoint operators. It is unital if it contains $\Id$. Let $\mathcal A$ be the smallest such unital subalgebra containing all $B_{uv}$.
\end{definition}
The preceding lemmas give
\begin{equation}\label{alg:Jsubset}
 \mathcal A\subset\{A=A^*:A\text{ is compatible with }K\}.
\end{equation}
Powers and real polynomials of any one member of $\mathcal A$ also belong
to $\mathcal A$: if $A^j\in\mathcal A$, then
$A^{j+1}=A\circ A^j$.

\begin{lemma}\label{alg:rankone}
The Jordan algebra $\mathcal A$ contains a rank-one orthogonal projection.
\end{lemma}
\begin{proof}
A nonzero symmetric tensor is determined by its diagonal polynomial:
the coefficient of $t_1t_2t_3t_4$ in
$B(\sum t_ix_i,\sum t_ix_i,\sum t_ix_i,\sum t_ix_i)$ is
$4!B(x_1,x_2,x_3,x_4)$.
Thus $B(v,v,v,v)$ is not identically zero.
After replacing $B$ by $-B$ if necessary, choose a unit vector $v$ where
this polynomial has a positive maximum $\beta$ on the unit sphere.
Write $A=B_{vv}$. The first and second variations along
\[
 v(t)=\cos(|w|t)v+\sin(|w|t)\frac{w}{|w|},
 \qquad w\perp v,
\]
give
\[
 Av=\beta v,\qquad
 0\geq-4\beta|w|^2+12h(Aw,w).
\]
Hence every eigenvalue on $v^\perp$ is at most $\beta/3<\beta$.
The eigenvalue $\beta$ is simple.
By the spectral theorem, a real interpolation polynomial in $A$ is exactly
the orthogonal projection $P_v$ onto $\R v$.
Since $A\in\mathcal A$ and polynomials in $A$ lie in $\mathcal A$,
we have $P_v\in\mathcal A$.
Changing the sign of $B$ leaves its generated Jordan algebra unchanged.
\end{proof}

\begin{lemma}\label{alg:rankone-generation}
An irreducible unital Jordan subalgebra of the self-adjoint operators that
contains a rank-one orthogonal projection is the entire space of self-adjoint
operators.
\end{lemma}
\begin{proof}
Write $P_v=v\otimes v^\flat$, with $|v|=1$, and define
\[
 H(u,w)=u\otimes w^\flat+w\otimes u^\flat,\qquad
 E=\{u:H(v,u)\in\mathcal A\}.
\]
The subspace $E$ contains $v$, since $H(v,v)=2P_v$.
For $A\in\mathcal A$,
\[
 2A\circ P_v=H(v,Av),
\]
so $Av\in E$.
If $u,w\in E\cap v^\perp$, direct multiplication gives
\[
 H(v,u)\circ H(v,w)
 =\tfrac12H(u,w)+h(u,w)P_v.
\]
Thus $H(u,w)\in\mathcal A$.
Subtracting components along $v$ shows this conclusion also holds for
arbitrary $u,w\in E$.

For $A\in\mathcal A$ and $u\in E$, the identity
\[
 2A\circ H(v,u)=H(Av,u)+H(v,Au)
\]
now shows $Au\in E$, because $Av,u\in E$.
Consequently $E$ is a nonzero invariant subspace for $\mathcal A$.
Irreducibility gives $E=H$.
The preceding product formula then yields every $H(u,w)$ in $\mathcal A$.
These operators span all self-adjoint operators.
\end{proof}

\begin{proof}[Proof of Theorem~\ref{algebra:isotropy}]
The hypothesis in Theorem~\ref{algebra:isotropy} that the self-adjoint common
commutant is trivial implies irreducibility of the contractions: a common invariant subspace
has an invariant orthogonal complement, because the contractions are
self-adjoint, and its orthogonal projection belongs to the self-adjoint
common commutant. Such a projection can only be $0$ or $\Id$.
The generated algebra $\mathcal A$ is therefore irreducible.
Lemmas~\ref{alg:rankone} and~\ref{alg:rankone-generation} imply
\[
 \mathcal A=\operatorname{Sym}(H,h).
\]
By \eqref{alg:Jsubset}, every self-adjoint operator is compatible with $K$.

For a unit vector $z$, compatibility of $P_z$ applied to $x,y\perp z$
gives
\[
 K(x,y)z=0.
\]
Choose an orthonormal basis $e_1,\ldots,e_m$.
For distinct indices $i,j,k$ this gives $K(e_i,e_j)e_k=0$.
Skew-adjointness then shows that $K(e_i,e_j)$ is supported on
$\operatorname{span}(e_i,e_j)$, so there are scalars $s_{ij}=s_{ji}$ such that
\[
 K(e_i,e_j)e_j=s_{ij}e_i,\qquad
 K(e_i,e_j)e_i=-s_{ij}e_j.
\]
If $m=2$, this already proves the result.
If $m\geq3$, apply the same observation to the mutually orthogonal vectors
$e_i$, $e_j+e_k$, $e_j-e_k$. It gives
\[
 0=K(e_i,e_j+e_k)(e_j-e_k)
   =(s_{ij}-s_{ik})e_i.
\]
Thus $s_{ij}=s_{ik}$. Symmetry in $i,j$ shows that all $s_{ij}$ equal one
scalar $s$. The displayed basis formulas determine $K$ completely and give
\[
 K(x,y)z=s\{h(y,z)x-h(x,z)y\}.
\]
This proves the theorem.
\end{proof}

\begin{lemma}\label{app:angular-equation}
For a smooth Landsberg metric, the angular restriction of the metric-lowered Berwald tensor satisfies~\eqref{eq:twistedCodazzi}.
\end{lemma}
\begin{proof}
Fix the base point and use the notation of~\eqref{eq:angularBerwaldDefinition}. By Lemma~\ref{lem:angularnonzero}, $B$ is fully symmetric and vanishes when the reference vector is inserted in any slot. All derivatives below are ordinary fibre derivatives in the original linear coordinates.

By full symmetry, the index $i$ may be taken as the metric-lowered index in $B_{lijk}$:
\[
B_{lijk}=g_{iq}\,\partial_{y^l}\partial_{y^j}\partial_{y^k}G^q.
\]
Commutativity of the four fibre differentiations cancels the fourth derivatives of $G$ in the following difference:
\begin{align}
\partial_pB_{lijk}-\partial_lB_{pijk}
&=2C_{piq}B^q{}_{ljk}
  -2C_{liq}B^q{}_{pjk}\notag\\
&=2\bigl(C_{pi}{}^qB_{lqjk}-C_{li}{}^qB_{pqjk}\bigr).
\label{canonical:ordinary}
\end{align}
Denote the expression in parentheses on the last line by $\Delta_i$. Symmetry of $B$ in $i,j,k$ makes the left side invariant under their permutations, so
\[
\Delta_i=\Delta_j=\Delta_k.
\]
The Levi--Civita connection of the fibre Hessian metric has coefficients $C_{ij}{}^k$. In the difference of its two covariant derivatives, the corrections on the $p,l$ slots cancel. The remaining three corrections are exactly $-\Delta_i-\Delta_j-\Delta_k$. Therefore
\begin{align*}
(\nabla^g_pB)_{lijk}-(\nabla^g_lB)_{pijk}
&=2\Delta_i-\Delta_i-\Delta_j-\Delta_k\\
&=-\Delta_i\\
&=C_{li}{}^qB_{pqjk}-C_{pi}{}^qB_{lqjk}.
\end{align*}
On the unit indicatrix the difference between $\nabla^g$ and $\nablah$ in angular directions is radial, and all radial insertions of $B$ vanish. Restricting to angular arguments therefore gives~\eqref{eq:twistedCodazzi}.

\end{proof}

\begin{proof}[Proof of Proposition~\ref{prop:angular-rigidity}]
Lemma~\ref{app:angular-equation} establishes the angular equation for the Berwald quartic. For an arbitrary nonzero quartic satisfying that equation, the compatibility identity~\eqref{eq:angularCompatibility} follows by symmetry. If $0<q<m$, Proposition~\ref{degenerate:conclusion} gives the rank restriction and the vanishing of $C$ on $U$.

If $q=m$, the full-rank set $\Omega$ is $U$ and is nonempty. Proposition~\ref{fullrank:generic} gives a point where the self-adjoint common commutant is $\R\Id$. Theorem~\ref{algebra:isotropy} then applies to $(T_p\Sigma,h_p)$ and yields $[C_X,C_Y]=sR_0(X,Y)$ at that point. This proves the proposition.
\end{proof}

\section*{Use of generative AI}
The authors used AI tools to assist in exploring proof strategies and checking calculations.
The authors are solely responsible for verifying all arguments and calculations and for the
final content.

\bigskip
\noindent\textsc{Zhuoqi Li}\\
Dongguan Middle School Songshan Lake School, Dongguan, Guangdong, China\\
\textit{E-mail address:} \href{mailto:429774494@qq.com}{429774494@qq.com}

\medskip
\noindent\textsc{Jianyu Mao}\\
School of Mathematical Sciences, Shanghai Key Laboratory of PMMP, East China Normal University, Shanghai 200241, China\\
\textit{E-mail address:} \href{mailto:1430658938@qq.com}{1430658938@qq.com}

\medskip
\noindent\textsc{Linfeng Zhou}\\
School of Mathematical Sciences, Shanghai Key Laboratory of PMMP, East China Normal University, Shanghai 200241, China\\
\textit{E-mail address:} \href{mailto:lfzhou@math.ecnu.edu.cn}{lfzhou@math.ecnu.edu.cn}

\end{document}